\documentclass[twoside,11pt]{article}

\usepackage{jmlr2e}
\usepackage{amsmath}
\allowdisplaybreaks[4]
\usepackage{amsbsy}
\usepackage{amscd}
\usepackage{amsopn}
\usepackage{amstext}
\usepackage{amsxtra}
\usepackage{epsfig}
\usepackage{times}

\def\lam{\lambda}

\jmlrheading{x}{20yy}{w-z}{a/bb}{cc/dd}{Balachandran, Kolaczyk, and Viles}

\ShortHeadings{On the Propagation of Low-Rate Measurement Error}{Balachandran, Kolaczyk, and Viles}
\firstpageno{1}

\begin{document}

\title{On the Propagation of Low-Rate Measurement Error to Subgraph Counts in Large Networks}

\author{\name Prakash Balachandran \email prakashb@bu.edu \\
              \name Eric D. Kolaczyk \email kolaczyk@bu.edu \\
              \addr Department of Mathematics \& Statistics, Boston University, Boston, MA 02215  USA 
              \AND
      \name Weston D. Viles \email weston.d.viles@dartmouth.edu \\
       \addr Department of Biomedical Data Science, Dartmouth College, Hanover, NH 03755 USA}
     
\editor{xxxxxxxxxxx}

\maketitle

\begin{abstract}
Our work in this paper is inspired by a statistical observation that is both elementary and broadly relevant to 
network analysis in practice -- that the uncertainty in approximating some true network graph $G=(V,E)$
by some estimated graph $\hat{G}=(V,\hat{E})$ manifests as errors in the status of (non)edges that must necessarily propagate
to any estimates of network summaries $\eta(G)$ we seek.  Motivated by the common practice of using plug-in estimates
$\eta(\hat{G})$ as proxies for $\eta(G)$, our focus is on the problem of characterizing the distribution of the
discrepancy $D=\eta(\hat{G}) - \eta(G)$, in the case where $\eta(\cdot)$ is a subgraph count.
Specifically, we study the fundamental case where the statistic of interest is $|E|$, the number of edges in $G$.
Our primary contribution in this paper is to show that in the empirically relevant setting of large graphs with 
low-rate measurement errors, the distribution of $D_E=|\hat{E}| - |E|$ is well-characterized by a Skellam distribution,
when the errors are independent or weakly dependent.  Under an assumption of independent errors, we are able to further show conditions under which this characterization is strictly better than that of an appropriate normal distribution.
These results derive from our formulation of a general result, quantifying the accuracy with which the difference of two sums
of dependent Bernoulli random variables may be approximated by the difference of two independent Poisson random 
variables, i.e., by a Skellam distribution.  This general result is developed through the use of Stein's method, and may be of 
some general interest.  We finish with a discussion of possible extension of our work to subgraph counts $\eta(G)$ of higher order.
\end{abstract}


\section{Introduction}
\label{intro}

The analysis of network data is widespread across the scientific disciplines.
Technological and infrastructure, social, biological, and information networks
are a few of the major network classes in which such analyses have been employed.  However, 
despite the already substantial body of work in network analysis generally (e.g., see~\citep{jackson2010social,kolaczyk2009statistical,newman2010networks}), 
with contributions from a variety of different fields of study, much work still
remains to more fully develop the theory and methods of statistical
analysis of network data, particularly for certain classes of problems of a fairly
fundamental nature. Here in this paper we pose and address a version of one such fundamental
problem, that regarding the propagation of error through the process of network
construction and summary.

In applied network analysis, a common {\it modus operandi} is to (i) gather basic measurements
relevant to the interactions among elements in a system of interest, (ii) construct a network-based
representation of that system, with nodes serving as elements and links indicating interactions between
pairs of elements, and (iii) summarize the structure of the resulting network graph using a variety of 
numerical and visual tools.  Key here is the point that the process of network analysis usually rests upon some
collection of measurements of a more basic nature.  For example, online social networks (e.g., Facebook) are based on the
extraction and merging of lists of `friends' from millions of individual accounts.  Similarly, biological networks (e.g., of gene
regulatory relationships) are often based on notions of association (e.g., correlation, partial correlation, etc.)
among experimental measurements of gene activity levels.  Finally, maps of the logical Internet traditionally have
been synthesized from the results of surveys in which paths along which information flows are learned
through a large set of packet probes (e.g., traceroute).

Importantly, while it is widely recognized that there is measurement error associated with these and other common
types of network constructions, most applied network analyses in practice effectively proceed as if there were in fact no error.
There are at least two possible reasons for this current state of affairs: (1) there is comparatively little
in the way of formal probabilistic analyses characterizing the propagation of such error and of statistical methods accounting 
for such propagation, and (2) in many settings (arguably due at least in part to (1)), much attention is given at the stages of
measurement and network construction to trying to keep the rate of error `low' in declaring the presence and absence
of links between nodes.

Here we offer a formal and general treatment of the problem 
of propagation of error, in which we provide a framework in which to characterize the manner in which 
errors made in assigning links between nodes accumulate in the reporting of certain functions of the network 
as a whole.  We provide a probabilistic treatment, wherein our goal is to understand the nature of the
distribution induced on the graph functions by that of the errors in the graph construction.

More formally, we consider a setting wherein an underlying (undirected)
network-graph $G$ possesses a network characteristic $\eta\left(G\right)$ of
interest.  While there are many types of functions $\eta(\cdot)$ used in practice to characterize networks
(e.g., centralities, path-based metrics, output from methods of community detection, etc.)
we restrict our attention here to the canonical problem of subgraph counting.  That is, we
are interested in the class of functions $\eta$ of the form
\begin{equation}
\eta_H(G) = \frac{1}{| Iso(H)|} \sum_{H'\subseteq K_{n_v},H'\cong H} 1{\{H'\subseteq G\}} \enskip ,
\label{eq:subg.cnt}
\end{equation}
where $n_v=|V(G)|$ is the number of vertices in $G$, $K_{n_v}$ is the complete graph on $n_v$ vertices,
$H$ is a graph of interest (i.e., copies of which we desire to count), and $H\subseteq G$ indicates
that $H$ is a subgraph of $G$ (i.e., $V(H)\subseteq V(G)$ and $E(H)\subseteq E(G)$).  The value
$|Iso(H)|$ is a normalization factor for the number of isomorphisms of $H$.  We will concentrate 
primarily on the fundamental case where $\eta(G)=|E|$, i.e., the number of edges in $G$.

If $\hat{G}$ is a network-graph resulting from an attempt to construct
$G$ from some collection of basic measurements, then the common practice of
reporting the analogous characteristics of $\hat{G}$ is equivalent to the use
of $\eta\left(\hat{G}\right)$, i.e. effectively a plug in estimator.  Let the discrepancy between
these two quantities be defined as $D=\eta\left(\hat{G}\right)-\eta\left(G\right)$, which in 
the case of counting edges reduces to $D_E=|\hat{E}|-|E|$.  Our 
goal is to make precise probabilistic statements about the behavior of $D$ under certain conditions.

Importantly, in the case where $\eta$ is defined as a subgraph count, as in (\ref{eq:subg.cnt}), $D$ may be 
expressed as the difference of (i) the number of times the subgraph $H$ arises somewhere in $\hat{G}$
but does {\em not} in fact exist in the same manner in $G$, and (ii) vice versa.  Hence, $D$
may be understood in this context to be the difference in total number of Type I and Type II
errors, respectively.  Intuitively, in the cases where a sufficiently low rate of error occurs on a large graph $G$, each of these
two sums should have a Poisson-like behavior.  This observation suggests that the propagation of low-rate measurement 
error to subgraph counts should behave, under appropriate conditions, like the difference of two Poisson random variables, i.e., 
a so-called Skellam distribution~\citep{skellam1946frequency}.  

Our contribution in this paper is to provide an initial set of results on the accuracy with which the Skellam distribution may be used in approximating the distribution of $D$, under the setting where the graph $G$ is  large and the rate of error is low.  We consider the cases of both sparse and dense networks. Our approach is through the use of Stein's method (e.g, \citep{barbour2005introduction}).  Specifically, we present a Stein operator for the Skellam probability distribution and, in a manner consistent with the Stein methodology, we derive several bounds on the discrepancy between the distribution of the difference of two arbitrary sums of binary random variables to an appropriately parameterized Skellam distribution. The latter in turn is then used to establish in particular the rate of weak convergence of $D_E$ to an appropriate Skellam random variable, under either independent or weakly dependent measurement errors.

As remarked above, there appears to be little in the way of a formal and general treatment of the error propagation
problem we consider here.  However, there are, of course, several areas in which the probabilistic or statistical treatment 
of uncertainty enters prominently in network analysis. The closest area to what we present here 
is the extensive literature on distributions of subgraph counts in random graphs.  See~\citep{janson2011random}, for example, 
for comprehensive coverage.  Importantly, there the graph $G$ is assumed to emerge according to a (classical) random graph
and uncertainty typically is large enough that normal limit theorems are the norm (although Poisson limit theorems
also have been established).  In contrast, in our setting we assume that $G$ is a fixed, true underlying graph, and then study the 
implications of observing a `noisy' version $\hat{G}$ of that graph, under various assumptions
on the nature of the noise, which involves two specific types of error (i.e., Type I and II errors), the contributions of which are informed 
in part by the topology of $G$ itself.  An area related to this work in random graphs is the work in statistical models for 
networks, such as exponential random graph models (ERGMs).  See~\citep{lusher2012exponential} for a recent 
treatment.  Here, while these models are inherently statistical in nature, the randomness due to generation of the
graph $G$ and due to observation of $G$ -- resulting in what we call $\hat{G}$ -- usually are combined into a single realization
from the underlying distribution.  And while subgraph counts do play a key role in traditional ERGMs, they typically enter
as covariates in these (auto)regressive models.  In a somewhat different direction, uncertainty in network construction 
due to sampling has also been studied in some depth.  See, for example, \citep[Ch 5]{kolaczyk2009statistical}
or~\citep{ahmed2012network} for surveys of this area.  However, in this setting, the uncertainty arises only from sampling 
-- the subset of vertices and edges obtained through sampling are typically assumed to be observed without error.  
Finally, we note that there just recently have started to emerge in the statistics literature formal treatments of the same type of graph observation error model that we propose here.  There the emphasis is on producing estimators of network model parameters and/or classifiers (e.g., \citep{priebe2012statistical}), for example, rather than on the type of basic network summary statistics that are the focus here.

The organization of this paper is as follows.  In Section~\ref{sec:initial.results} we provide necessary background.
In Section~\ref{sec:general.results} we then provide a general set of results useful
for our general problem.  Specifically, we establish a bound for the Kolmogorov-Smirnov distance 
between the distribution of the difference of two arbitrary sums of binary random variables from a certain Skellam.  
This work is based on the application of Stein's method to the Skellam distribution, a first of its kind to the
best of our knowledge, and the results therefore are of some independent interest as well.  
In Section~\ref{sec:applications} we then illustrate the way in which these general results
may be used to understand the propagation of error in networks for counting edges.
In doing so, several other general results are established.  Some implications of these results on the 
problem of counting subgraphs of higher order are noted in Section~\ref{sec:discussion}, along with 
other related discussion.   Proofs of our key results may be found in the appendices.  

\section{Background}
\label{sec:initial.results}

\subsection{Notation and Assumptions}

By $G=(V,E)$ we will mean an undirected graph, with vertex set $V$ of cardinality $|V|$ and edge set $E$ of cardinality $|E|$.  Much of the results that follow will be stated as a function of the number of vertices which, for notational convenience, we denote $n_v=|V|$.  Let $\mu=\mu(G) = 2|E|/n_v$ correspond to the average degree of a vertex in $G$.  We assume the vertex set $V$ is known but that the edge set $E$ is unknown.  However, we assume there is information by which to construct an approximation to $E$ or, more formally, by which to infer $E$, as a set $\hat{E}$, yielding an inferred network graph $\hat{G}=(V,\hat{E})$.  

While there are many ways in practice by which the set $\hat{E}$ is obtained, one principled
way of viewing the task is as one of performing 
$\binom{n_v}{2}$ hypothesis tests, using the data underlying the graph construction process as input,
one for each of the vertex pairs $\left\{i,j\right\}$ in
the network graph $G$.  In some contexts, $\hat{G}$ is literally obtained through hypothesis
testing procedures; for instance, in constructing some gene regulatory networks from microarray
expression data.  See~\citep[Ch 7]{kolaczyk2009statistical}, for example.
Formally, in such cases we can think of $\hat{G}$ as resulting from a collection of testing problems
\begin{equation}\nonumber
H_0:\left\{i,j\right\}\notin E\mbox{ versus }H_1:\left\{i,j\right\}\in E \enskip ,
\end{equation}
for $\{i,j\}\in V^{(2)}$, where
\begin{equation}\nonumber
V^{\left(2\right)}=\left\{\left\{i,j\right\}:i,j\in V;i<j\right\} \enskip .
\end{equation}
These tests amount to a collection of $\binom{n_v}{2}$ binary random variables $\left\{Y_{ij}:\left\{i,j\right\}\in V^{\left(2\right)}\right\}$, where
\begin{equation}\nonumber
Y_{ij}=\left\{\begin{array}{ll} 1 & \mbox{if $H_0$ is rejected}\\ 0 & \mbox{if $H_0$ is not rejected.}\end{array}\right.
\end{equation}
Note that the random variables $Y_{ij}$ need not be independent and, in fact, in many contexts will most likely be dependent.  Gene regulatory networks inferred by correlating gene expression values at each vertex $i$ with that of all other vertices $j\in V\setminus \{i\}$ and maps of the logical Internet obtained through merging paths learned by sending traffic probes 
between many sources and destinations are just two examples where dependency can be expected.

Whether obtained informally or formally, we can define the inferred edge set $\hat{E}$ as
\begin{equation}\nonumber
\hat{E}=\left\{\left\{i,j\right\}\in V^{\left(2\right)}:Y_{ij}=1\right\} \enskip .
\end{equation}
It is useful to think of the collection of random variables
$\left\{Y_{ij}:\left\{i,j\right\}\in V^{\left(2\right)}\right\}$ 
as being associated with two types of errors.  That is, we express the marginal
distributions of the $Y_{ij}$ in the form
\begin{equation}
  Y_{ij}\, \sim \, \begin{cases}
  \hbox{Bernoulli}\left(\alpha_{ij}\right) , \text{if $\{i,j\}\in E^c$,}\\
  \hbox{Bernoulli}\left(1 - \beta_{ij}\right), \text{if $\{i,j\}\in E$,}
\end{cases}
\label{eq:basic.bernoulli}
\end{equation}
where $E^c = V^{(2)}\setminus E$.  
Again pursuing the example of network construction
based on hypothesis testing, $\alpha_{ij}$ can be interpreted as 
the probability of Type-I error for the test on vertex pair
$\left\{i,j\right\}\in E^c$, while $\beta_{ij}$
is interpreted as the probability of Type-II error for the test on vertex pair
$\left\{i,j\right\}\in E$.

Our interest in this paper is in characterizing the manner in which the uncertainty in the $Y_{ij}$ propagates to 
subgraph counts on $\hat{G}$.  More specifically, we seek to characterize the distribution of the difference
\begin{equation}
D = \frac{1}{| Iso(H)|} \sum_{H'\subseteq K_{n_v},H'\cong H} \left[ 1{\{H'\subseteq \hat{G}\}} - 1{\{H'\subseteq G\}}\right] \enskip ,
\label{eq:subg.cnt.diff}
\end{equation}
for a given choice of subgraph $H$.  Naturally, this distribution will depend in no small part on context.  Here we focus on
a general formulation of the problem in which we make the following three assumptions.
\begin{itemize}
  \item[({\bf A1})] \quad {\em Large Graphs:}  $n_v \rightarrow \infty$.
  \smallskip
  
  \item[({\bf A2})]  \quad {\em Edge Unbiasedness:} $\sum_{\{i,j\}\in E^c} \alpha_{ij} = \sum_{\{i,j\}\in E} \beta_{ij} \left( \equiv \lambda\right)$.
   \smallskip

  \item[({\bf A3})]  \quad {\em Low Error Rate:} $\lambda = \Theta\left( \mu \right)$.
\end{itemize}
\medskip

Assumption (A1) reflects both the fact that the study of large graphs is a hallmark of modern applied work in complex networks and, accordingly, our desire to make statements that are asymptotic in $n_v$.  

Our use of assumption (A2) reflects the understanding that a `good' approximation $\hat{G}$ to the graph $G$ should at the very least have roughly the right number of edges.  The difference of the two sums defined in (A2) is in fact the expectation of the statistic $D$ in (\ref{eq:subg.cnt.diff}) for the case where the subgraph being counted is just a single edge, i.e., it is the expected discrepancy between the number of observed edges $|\hat{E}|$ and the actual number of edges $|E|$.  So (A2) states that this particular choice of $D$ have expectation zero.  Alternatively, (A2) may be interpreted as saying that the total numbers of Type I and Type II errors should be equal to a common value $\lambda$.  

Finally, in (A3) we encode the notion of there being a `low' rate of error in $\hat{G}$.  Specifically, we assume that the number of Type I or Type II errors in edge status that we expect throughout the network is roughly on par with the average number of edges incident to any given vertex in the network.  This condition can be re-expressed in a useful manner with respect to $n_v$ if, as is common in the literature, we distinguish between  sparse and dense graphs.  By the term {\em sparse} we will mean a graph for which $|E|= \Theta\left(n_v \log n_v\right)$, and by {\em dense}, $|E|=\Theta\left( n_v^2\right)$.  Hence, assumption (A3) reduces to $\lambda=\Theta\left(\log n_v\right)$ in the case of sparse graphs, and to $\lambda=\Theta(n_v)$, in the case of dense graphs.

In addition, for convenience, we add to the core assumptions (A1)-(A3) a fourth assumption, upon which we will call periodically throughout the paper when desiring to simplify some of our expressions.
\begin{itemize}
  \item[({\bf A4})] \quad {\em Homogeneity:} $\alpha_{ij}\equiv \alpha$ and $\beta_{ij}\equiv \beta$, for $\alpha,\beta\in (0,1)$.
\end{itemize}
\medskip
In other words, we assume that the probability of making a Type I or II error (as the case may be) does not depend upon
the specific (non)edge in question.

Lastly, for completeness, we recall the definition of the Skellam distribution.
A random variable $W$ defined on the integers is said to have a Skellam distribution, with parameters $\lambda_1, \lambda_2 > 0$, i.e.,
$W\sim\mbox{Skellam}\left(\lambda_1,\lambda_2\right)$, if 
\begin{equation} 
\mathbb{P}\left(W=k\right)=e^{-\left(\lambda_1+\lambda_2\right)}\left(\sqrt{\frac{\lambda_1}{\lambda_2}}\right)^kI_k\left(2\sqrt{\lambda_1\lambda_2}\right)\mbox{
for $k\in\mathbb{Z}$},
\label{eq:skellam}
\end{equation}
where $I_k\left(2\sqrt{\lambda_1\lambda_2}\right)$ is the modified Bessel
function of the first kind with index $k$ and argument
$2\sqrt{\lambda_1\lambda_2}$.  The Skellam distribution may be constructed 
by defining $W$ through the difference of two independent Poisson random variables, with means $\lambda_1$ and $\lambda_2$, 
respectively.  The mean and variance of this distribution are given by $\mathbb{E}[W] = \lambda_1-\lambda_2$
and $\hbox{Var}(W) = \lambda_1 + \lambda_2$.  The distribution of $W$ is in general nonsymmetric, with symmetry
holding if and only if $\lambda_1 = \lambda_2$.

The main results we provide in this paper are in the form of bounds on the extent to which the distribution of 
random variables like the discrepancy $D$ in (\ref{eq:subg.cnt.diff}) may be well-approximated by an appropriate Skellam distribution.  For this purpose, we adopt the Kolmogorov-Smirnov distance to quantify the distance between distributions of two random variables, say, $X_1$ and $X_2$, i.e., 
\begin{equation}\label{ks}
ds_{KS}(X_1,X_2) \equiv \sup_{x\in\mathbb{R}}\left|\mathbb{P}\left(X_1\leq x\right)-\mathbb{P}\left(X_2\leq x\right)\right|.
\end{equation}

\subsection{Counting Edges}
\label{sec:motiv.examples}

Generic subgraph counts, and the corresponding noise in obtaining them, can be quite varied in real applications.  
Accordingly, most of our specific results pertain to the fundamental case of counting edges.  That is, where
the choice of subgraph $H$ is simply a single edge, and therefore the function $\eta(G)$ in (\ref{eq:subg.cnt}) is just 
the total number of edges in $G$, i.e., $\eta(G)=|E|$.  We will consider two scenarios for this case, wherein the random variables $Y_{ij}$ are independent or weakly dependent.  

In the case where the edge noise is independent, the discrepancy 
\begin{eqnarray}
D_E & = & |\hat{E}| - |E| \nonumber \\
      & = & \sum_{\left\{i,j\right\}\in E^c}Y_{ij}-\sum_{\left\{i,j\right\}\in E}\left(1-Y_{ij}\right) \enskip 
\label{eq:diff.edge.cnt}
\end{eqnarray}
has expectation $\mathbb{E}[D_E] = \alpha |E^c| - \beta |E| = \lam - \lam = 0$ and
variance $\sigma^2 = \alpha(1-\alpha)|E^c| + \beta(1-\beta)|E|$, and its behavior can be established
using existing methods from the literature (i.e., essentially, Chen-Stein methods).  However, we include it as
an important base case, comparing results obtainable by our methods to those obtainable by more traditional
techniques, in Section~\ref{sec:edge.cnts.w.indep}.

Alternatively, suppose that the variables $Y_{ij}$ are {\em dependent}.  The random variable $D_E$ again has expectation zero, although its variance -- and hence its asymptotic behavior -- will differ from the independent case,
in a manner dictated by the nature of the underlying dependency in the noise.  It often can be expected in practice that the error associated with construction of the empirical graph $\hat{G}$ will involve dependency across (non)edges.  For
example, relations in gene regulatory networks are often declared based on sufficiently strong levels of association 
between gene-specific measurements (e.g., measures of gene expression).  The comparison of the measurements
for each gene with those of all of the others necessarily induces potential dependencies among the random variables $Y_{ij}$. However, a precise characterization 
of such dependency is typically problem-specific and, more often than not, nontrivial in nature.  
In Section~\ref{sec:edge.cnts.w.dep} we will assume general dependency conditions in the spirit of 
traditional monotone coupling arguments, which will allow for further analysis and interpretation.

\section{General Results on Approximation by Skellam}
\label{sec:general.results}

Recall the general form of our statistic of interest $D$ in (\ref{eq:subg.cnt.diff}), as the difference
of two sums of binary random variables.  Under appropriate conditions it seems reasonable to expect that 
the distribution of $D$ be well-approximated by a Skellam distribution.  And for the simplest case, in which
we are counting edges under independent noise, it is possible to show that this is in fact the case,
through manipulation of existing results for approximating sums by Poisson distributions.  Without independence, however, it is necessary to approach the problem directly, by explicitly using the Skellam distribution.  In this
section, we therefore provide the results of such an approach.  This is a completely general treatment
-- devoid of the motivating context of counting subgraphs -- and therefore also likely of some independent interest. 
In Section~\ref{sec:applications} we return to the problem of counting subgraphs under low-rate error
and illustrate the use of the results presented here in this section through application to the case of counting edges.

Our approach in this section is through Stein's method.  This choice is reminiscent, naturally, of the Chen-Stein treatment for Poisson approximations.  However, the task is technically more involved at several points, as it requires handling a Stein function that is defined through a second-order difference, rather than the first-order difference encountered in the Poisson problem.  Moreover, the kernel of the Skellam distribution includes a modified Bessel function of the first kind, which emerges in ways necessitating a somewhat delicate treatment.

\subsection{A Stein Bound for the Skellam Distribution}
\label{sec:stein.bnd}

Let $U$ be a random variable defined as 
\begin{equation} \label{V Def}
U = \sum_{k=1}^nL_k-\sum_{k=1}^mM_k \enskip ,
\end{equation}
where $\left\{\left\{L_k\right\}_{k=1}^n,\left\{M_k\right\}_{k=1}^m\right\}$ is a
collection of two sets of indicator random variables with $\mathbb{E}[L_k]=p_k$ for
$k=1,\ldots,n$ and $\mathbb{E}[M_k]=q_k$ for $k=1,\ldots,m$.  
In the case of our subgraph counting problem, $U=D$, where $D$ is defined in (\ref{eq:subg.cnt.diff}),
although for the remainder of this current section $U$ is defined generally.

Recall the definition of a Skellam random variable $W$ in (\ref{eq:skellam}).  We desire a bound on
\begin{equation}\label{ks}
d_{KS}(U,W):=\sup_{x\in\mathbb{R}}\left|\mathbb{P}\left(U\leq x\right)-\mathbb{P}\left(W\leq x\right)\right| \enskip ,
\end{equation}
quantifying how close the distribution of $U$ is to that of $W$.  
In pursuing the standard paradigm for Stein's method, we first determine an operator
$\mathcal{A}\left[f\left(k\right)\right]$ such that
\begin{equation} \nonumber
\begin{array}{ccc}
\mathbb{E}\mathcal{A}\left[f\left(W\right)\right]=0 & \mbox{ if and only if } & W\sim\mbox{Skellam}\left(\lambda_1,\lambda_2\right)
\end{array}
\end{equation}
for any bounded function $f:\mathbb{Z}\mapsto\mathbb{R}$.  This operator need not be unique, but the theory only requires one.  This is
accomplished through the following result, the proof of which uses several properties of the modified Bessel function of the first kind,
as detailed in the appendix, in Section~\ref{sec:proof.of.steinop}.
\begin{theorem}\label{thm:steinop}
A Stein operator $\mathcal{A}$ for the
$\mbox{Skellam}\left(\lambda_1,\lambda_2\right)$ distribution is
\begin{equation}\nonumber
\mathcal{A}\left[f\left(k\right)\right]=\lambda_1f\left(k+1\right)-kf\left(k\right)-\lambda_2f\left(k-1\right).
\end{equation}
\end{theorem}

With this operator in hand, and again following the usual paradigm under Stein's method, we set
\begin{equation} \nonumber
\mathcal{A}\left[f\left(k\right)\right]=g\left(k\right)
\end{equation}
for a class of test functions $g\left(k\right)$, and allow that to implicitly define the function $f$.  The choice of the test functions $g$ is guided by the choice of the metric used to measure the distance between $U$ and $W$.  Since the metric we choose to measure the distance between $U$ and $W$ is given by $d_{KS}(U,W)$ in (\ref{ks}), we choose the test function $g:=g_x$ given by
\begin{equation} \label{test function}
g_x\left(k\right) = 1\left\{k\leq x\right\}-\mathbb{P}\left(W\leq x\right)
\end{equation}
for any $x\in \mathbb{R}$.  

At this point it is common to exhibit a solution $f$ defined by our choice of $g$.  Instead, we forestall
that step until later in this section, choosing rather to state a general result that will allow us to more quickly
gain insight into the nature of the bounds we are able to obtain.  Our result employs a minor variant of the notion of coupling
that is common to the literature on Chen-Stein approximations.
\begin{theorem}\label{thm:coupling}
Let $U$ be as in (\ref{V Def}) and let $\mathcal{L}(U)$ denote the law of $U$.  Let $$\mathcal{L}\left(U_k^{(L)}+1\right)=\mathcal{L}\left(U | L_k=1\right) \; \; \; {\rm for} \; \; \; k=1,\ldots,n$$ and $$\mathcal{L}\left(U_k^{(M)}-1\right)=\mathcal{L}\left(U | M_k=1\right), \; \; \; {\rm for} \; \; \; k=1\ldots,m$$ be a collection of random variables all defined on a common probability space.  If $\lambda_1=\sum_{k=1}^n p_k$ and $\lambda_2=\sum_{k=1}^m q_k$, and $W\sim Skellam(\lambda_1,\lambda_2)$, then
\begin{equation}\label{Stein KS Bound}
d_{KS}(U,W)\leq ||\Delta f|| \left\{\sum_{k=1}^np_k \mathbb{E}\left|U-U_k^{(L)} \right|+\sum_{k=1}^mq_k\mathbb{E}\left| U-U_k^{(M)}
 \right| \right\}  \enskip ,
\end{equation}
where
$$||\Delta f|| =\sup_{x\in \mathbb{R}} \sup_{j\in\mathbb{Z}}\left|f_x\left(j+1\right)-f_x\left(j\right)\right|$$
and $f_x$ is a solution to $\mathcal{A}\left[f_x(k)\right]=g_x(k)$ for $k\in \mathbb{Z}$.
\end{theorem}

The proof of this result relies on elementary considerations of the equation $\mathcal{A}[f_x(k)]=g_x(k)$
and may be found in the appendix in Section~\ref{sec:proof.of.coupling}.
The extent to which it allows one to obtain error estimates of practical interest in a particular setting 
will depend on the extent to which both the main expression within brackets in (\ref{Stein KS Bound})
and the quantity $||\Delta f||$ can be further controlled. While control of the former is problem dependent, 
control of the latter is not, and may be dealt with separately, as we do next.  Afterwards, in 
Section~\ref{sec:applications}, we illustrate the control of the bracketed expression in 
(\ref{Stein KS Bound}), in the context of the problem of counting edges introduced in Section~\ref{sec:motiv.examples}.

\subsection{Controlling the term $||\Delta f||$.}
\label{sec:bounding.constant}

Controlling  $||\Delta f||$ in (\ref{Stein KS Bound}) first requires
understanding the solution $f_x(k)$.  We provide a family of closed-form expressions for this solution 
in the following.
\begin{theorem}\label{steinsol}
Let $g_k$ be defined as in equation (\ref{test function}).  If $f_x$ is a bounded solution to the difference equation
\begin{equation} \nonumber
\lambda_1f_x\left(k+1\right)-kf_x\left(k\right)-\lambda_2f_x\left(k-1\right)=g_x\left(k\right)
\end{equation}
for $k\in \mathbb{Z}$, then $f_x$ is given by
\begin{equation} \nonumber
f_x\left(m\right) = \left\{\begin{array}{ll}
\left(-1\right)^m\left(\sqrt{\frac{\lambda_2}{\lambda_1}}\right)^mI_m\left[\left(-1\right)^c\left(\sqrt{\frac{\lambda_1}{\lambda_2}}\right)^c\frac{1}{I_c}f\left(c\right)\right. & \\
\left.+\frac{e^{\lambda_1+\lambda_2}}{\sqrt{\lambda_1\lambda_2}}\sum_{n=c}^{m-1}\frac{\left(-1\right)^{n+1}}{I_nI_{n+1}}\mathbb{P}\left(W\leq\min\left\{n,x\right\}\right)\mathbb{P}\left(W>\max\left\{n,x\right\}\right)\right] & \mbox{if $m>c$}\\
\left(-1\right)^m\left(\sqrt{\frac{\lambda_2}{\lambda_1}}\right)^mI_m\left[\left(-1\right)^c\left(\sqrt{\frac{\lambda_1}{\lambda_2}}\right)^c\frac{1}{I_c}f\left(c\right)\right. & \\
\left.-\frac{e^{\lambda_1+\lambda_2}}{\sqrt{\lambda_1\lambda_2}}\sum_{n=m}^{c-1}\frac{\left(-1\right)^{n+1}}{I_nI_{n+1}}\mathbb{P}\left(W\leq\min\left\{n,x\right\}\right)\mathbb{P}\left(W>\max\left\{n,x\right\}\right)\right] & \mbox{if $m<c$}.
\end{array}\right.
\end{equation}
for any initial condition $\left(c,f_x\left(c\right)\right)$ with
$c\in\mathbb{Z}$ and $f_x\left(c\right)\in\mathbb{R}$.
\end{theorem}
\noindent The proof of this theorem is similar to that of solving a second order linear differential equation.  An integrating factor is found, integration is performed with a boundary condition at $-\infty$, and then a second integration is performed with the initial condition.
Details are provided in the appendix in Section~\ref{sec:proof.of.steinsol}.

Leveraging our insight into $f_x$ to control $||\Delta f||$ means producing 
a bound on the absolute differences $|\Delta f_x(j)|=|f_x(j+1)-f_x(j)|$ independent of $x\in \mathbb{R}$ and $j\in \mathbb{Z}$.  Consider first the special case where $\lambda_1$ and $\lambda_2$ are equal, 
for which we are able to offer the following result.
\begin{theorem}\label{lipschitz}
Suppose that $\lambda_1=\lambda_2 \equiv \lambda$, so that $\mathbb{E}[U]=0$
and $W$ is a $\hbox{Skellam}(\lambda,\lambda)$ distribution, symmetric about zero.  
Assume $ \lambda \ge 1$.  Then
\begin{equation} \nonumber
||\Delta f|| \leq  \frac{160}{2\lambda}.
\end{equation}
\end{theorem}
\noindent The proof of this theorem is highly technical in nature, and relies on a concentration inequality for the Skellam$(\lambda,\lambda)$ distribution \citep{kash} with several other technical arguments. A sketch of the proof may be found in the appendix in 
Section~\ref{sec:proof.of.lipschitz}, while a detailed presentation is available in the Supplementary Materials.

Note that the bound in Theorem~\ref{lipschitz} is essentially the analogue of the classical result for Poisson approximation,
in which, for sufficiently large $\lambda$, the term $1/\lambda$ is the standard factor.  In both cases, therefore, the
corresponding term $||\Delta f||$ is bounded by the inverse of the expected total number of counts, where here 
that is $\mathbb{E}[T_1+T_2] = 2\lambda$.

The above result is of immediate relevance to the problem of counting edges, given assumption (A2), whether under
the assumption that the edge noise is independent or dependent.  We will make use of this result in the next section.  For applications involving higher-order subgraphs, we can expect to need an extension of Theorem~\ref{lipschitz} to the general case of $\lambda_1\ne \lambda_2$.  For arbitrary $\lambda_1, \lambda_2>0$, we are unable to produce a satisfactory bound.  However, supported by preliminary numerical work,  we have the following conjecture.
\begin{conjecture}\label{lipschitz:l1nel2}
In general, for $\lambda_1$ and $\lambda_2$ sufficiently close and large,
\begin{equation} \nonumber
||\Delta f|| \leq  \frac{C}{\lambda_1+\lambda_2},
\end{equation}
for some constant $C>0$ independent of $\lambda_1, \lambda_2$.
\end{conjecture}
\noindent

\section{Application of General Results to Counting Edges}
\label{sec:applications}

We now illustrate the use of our general results for the problem of characterizing the propagation of low-rate measurement error to subgraph counts in large graphs, for the specific case of counting edges.  

\subsection{Edge Counts Under Independent Edge Noise}
\label{sec:edge.cnts.w.indep}

Recall the problem wherein the function of interest (\ref{eq:subg.cnt}) counts the number of edges 
in $G$, i.e., $\eta(G)=|E|$,  and the variables $Y_{ij}$ in (\ref{eq:basic.bernoulli}) are independent.  
In light of Theorems~\ref{thm:coupling} and~\ref{lipschitz}, we have the
following result characterizing the behavior of the discrepency $D$ in (\ref{eq:subg.cnt.diff}),
which here is simply $D_E = |\hat{E}| - |E|$.  
\begin{theorem}
Under assumptions (A1)-(A4) and independence among errors in declaring (non)edge status (i.e., among the $Y_{ij}$),
\begin{equation}
d_{KS}\left(\, D_E\, ,\,  \hbox{Skellam}(\lambda,\lambda)\, \right) \le O\left( n_v^{-1}\right)\enskip ,
\label{eq:skellam.error.indep}
\end{equation}
for both sparse and dense graphs $G$.
\label{thm:skellam.indep.case}
\end{theorem}

Proof of this result may be found in the appendix, in Section~\ref{sec:proof.of.normal.vs.skellam.indep}.
The theorem establishes a rate at which -- in large networks, whether sparse or dense, with independent and homogeneous low-rate errors -- the distribution of the discrepancy $D_E$ tends to that of an appropriate Skellam distribution, i.e.,  symmetric and centered on zero, with variance $2\lambda$.  The same rate may be established using more standard arguments from Chen-Stein theory, the proof of which we also include in the appendix, for completeness.  These latter arguments, of course, only hold in the case of independence assumed here, and do not extend generally to the case of dependence in the edge noise.

To put the rate established in the above theorem in better context, it is interesting to compare to the case where a normal distribution is used instead to approximate that of the discrepancy $D_E$.  The following theorem, proof of which also may be found in Section~\ref{sec:proof.of.normal.vs.skellam.indep}, provides both upper and lower bounds.
\begin{theorem}
Let $\sigma^2 = \hbox{Var}(D_E)$.  Under the same conditions as Theorem~\ref{thm:skellam.indep.case} , in the case of sparse graphs
\begin{equation}
\Omega\left( \log^{-1}n_v\right) \le  d_{KS}\left(\, D_E/\sigma \, , \, N(0,1)\, )\right) 
                         \le O\left( \log^{-1/2} n_v\right) \enskip ,
\enskip 
\label{eq:normal.error.indep.sparse}
\end{equation}
while in the case of dense graphs, 
\begin{equation}
\Omega\left(n_v^{-1}\right) \le d_{KS}\left(\, D_E/\sigma \, , \, N(0,1)\, )\right) 
          \le \Omega\left( n_v^{-1/2}\right)  \enskip ,
\label{eq:normal.error.indep.dense}
\end{equation}
where $N(0,1)$ refers to a standard normal random variable.
\label{thm:normal.indep.case}
\end{theorem}

\begin{figure}[bht]
  \centerline{\includegraphics[width=2.5in]{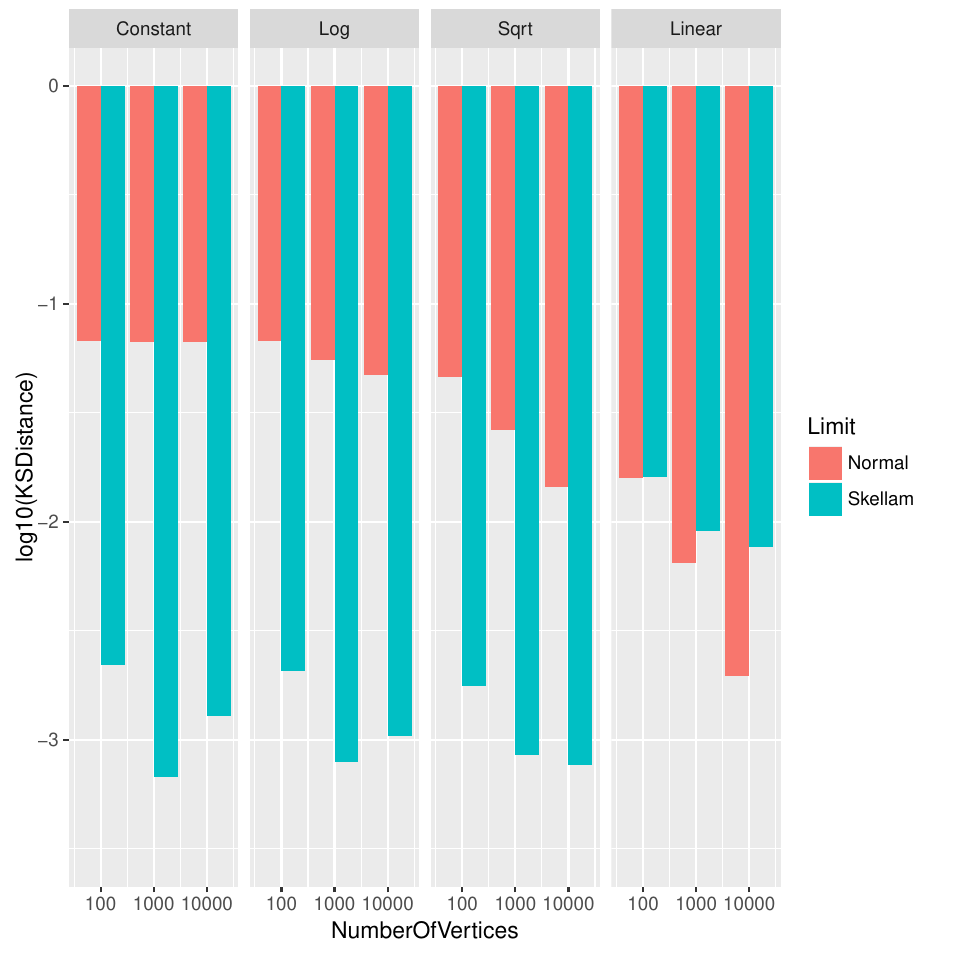}\quad
                  \includegraphics[width=2.5in]{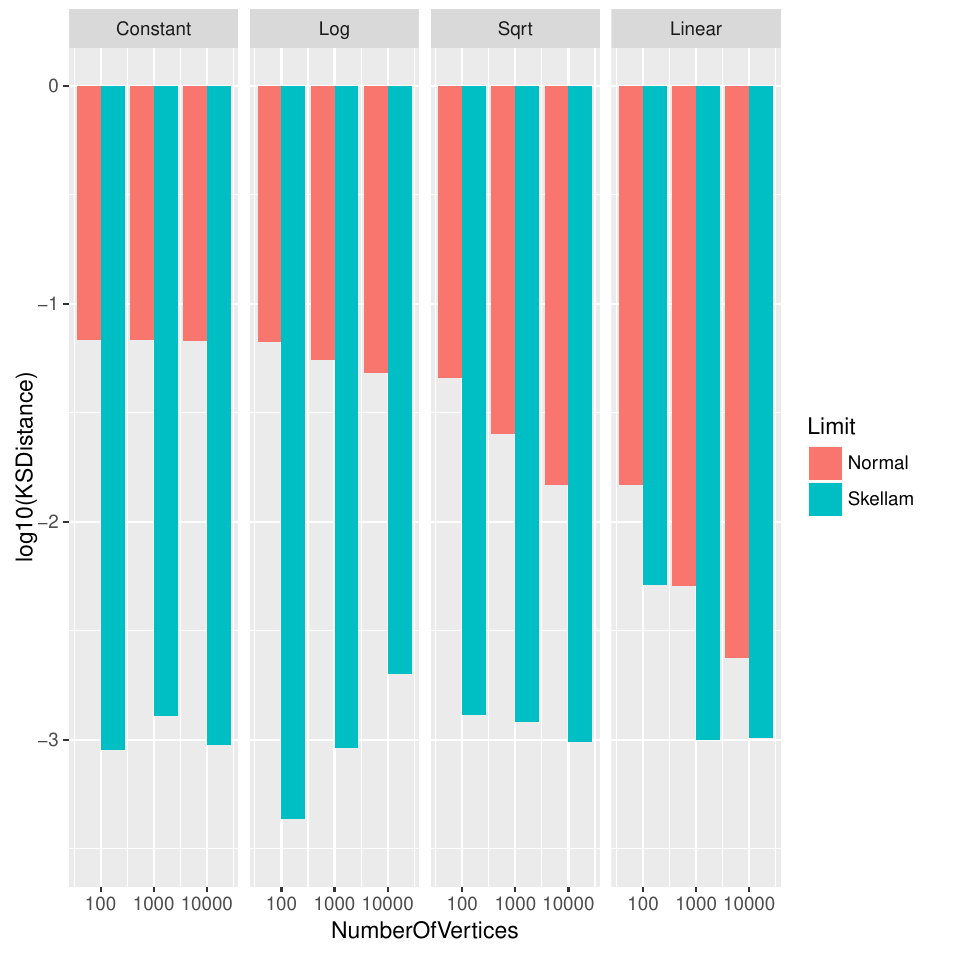}}
\caption{(Log)Kolmogorov-Smirnov distance between Skellam and standard normal approximations to the distributon of discrepency $D_E$ in edge counts under independent errors. {\em Left:} Sparse case. {\em Right:} Dense case.}
\label{fig:iid.simul}
\end{figure}

These two theorems together indicate that in this context a Skellam approximation is clearly superior to a normal for sparse graphs, and they suggest that it can be better as well for dense graphs.  These statements are supported by the results of a small simulation study, shown in Figure~\ref{fig:iid.simul}.  There we compared the two approximations as $n_v$ ranges from $100$ to $1000$ to $10,000$, for error rates $\lambda$ defined to be constant, logarithmic, square root, or linear functions of $n_v$.  For the sparse and dense cases, we let $|E|$ equal $n_v\log n_v$ and $n_v(n_v-1)/4$, respectively.  Looking at the sparse case, for when $\lambda=\log n_v$,  the Skellam approximation clearly dominates the normal.  However, interestingly, this dominance continues even when the error rate is set equal to $n_v^{1/2}$.  Only once the error rate is $n_v$ do we see the normal approximation begin to overtake the Skellam approximation. Note that by this stage, $\beta = O(1)$, and so essentially there is no `signal' standing out from the `noise'.  Similarly, looking at the dense case, we see that the Skellam approximation is better than the normal approximation at all error rates, including, in particular, when the error rate equals $n_v$, the case addressed by the above two theorems.

In summary, in the independent case, the Skellam distribution dominates the normal as an approximation when there can be expected to be a clear graph `signal' standing out against the `noise' induced by underlying low-rate measurement errors.  

\subsection{Edge Counts Under Dependent Edge Noise}
\label{sec:edge.cnts.w.dep}

Again, as just above, consider the context wherein counting edges is of interest, so that
$\eta(G)=|E|$ and our goal is to characterize the accuracy with which $D_E=|\hat{E}| - |E|$ is 
approximated by a $\hbox{Skellam}(\lambda,\lambda)$ random variable.  Now, however, we assume
that the error associated with construction of the empirical graph $\hat{G}$ will involve dependency across (non)edges. 
That is, the random variables $Y_{ij}$ are now dependent. 
A precise characterization of such dependency is typically problem-specific and, more often than not, nontrivial in nature.  
Here, for the purposes of illustration, we instead provide certain results of a general nature, working from the 
bound (\ref{Stein KS Bound}) of Theorem~\ref{thm:coupling}.

Of the two terms in (\ref{Stein KS Bound}), the first term $||\Delta f||$ is again known to behave as $O(\lambda^{-1})$, 
by Theorem~\ref{lipschitz}.  On the other hand, control of the second term, in brackets, requires some care.
For example, naive inter-change of absolute values and summations with expectation yields that
$$p_k \mathbb{E}\left|U-U_k^{(L)} \right| \le
p_k^2 + \sum_{k\ne j} p_kp_j + E[L_k L_j]  + \sum_{\ell =1}^m p_k q_\ell + E[L_k M_\ell] \enskip ,
$$
and similarly for $q_k\mathbb{E}\left| U-U_k^{(M)} \right|$.  Unfortunately, it is straightforward to show that
for the dependent error version of the problem considered in the previous subsection (i.e., involving 
independent and homogeneous low-rate errors on large-spare networks) the bound we obtain for 
$d_{KS}\left(\, D_E\, ,\,  \hbox{Skellam}(\lambda,\lambda)\, \right)$ will be
no better than $O(\lambda)$ -- regardless of the nature of the dependency among the $Y_{ij}$.

One possible approach to a more subtle handling of these terms is motivated by considerations of 
hypothesis testing.  Suppose that the $L_k$ correspond to indicators of Type I error for $n$ tests
under their corresponding null hypotheses, and the $M_k$, to indicators of Type II error for $m$ tests
under their corresponding alternative hypotheses.  Furthermore, suppose that the corresponding test statistics 
are all defined on the same scale and compared to the same threshold.  Moreover, for simplicity, we assume these statistics
all have non-negative values and that their distribution under the null sits to the left of that under the alternative, so that 
more extreme positive values tend to support the alternative.  In this setting, if we know, for example, 
that $L_1=1$, we know that at least one rejection of a null hypothesis has occurred, indicating that the threshold
sits to the left of the right-most extreme of the empirical null distribution.  Accordingly, we are inclined to 
expect that there may be other such rejections of the null, i.e., other Type I errors.  At the same time, we would expect fewer Type II errors, i.e., fewer $M$ that equal $1$.  Conversely, if we see a Type II error, 
say $M_1=1$, it can be argued that we would be inclined to expect more Type II errors and, at the same time, fewer Type I errors.

Together these high-level arguments suggest that a reasonable generic model for these errors
is one in which there is positive correlation {\em within} the $L$'s and $M$'s, respectively, but
negative correlation {\em between}.  The conditions of the following theorem capture this notion,
which in turn allow us to produce a sensible bound, improving on that of Theorem~\ref{Stein KS Bound}.
\begin{theorem}
Let $\tilde{L}^{L_k}_j$ and $\tilde{M}^{L_k}_\ell$ be random variables distributed as $L_j$ and $M_\ell$
respectively, conditional on $L_k=1$.  Similarly, let $\tilde{L}^{M_k}_j$ and $\tilde{M}^{M_k}_\ell$ be
distributed as $L_j$ and $M_\ell$, conditional on $M_k=1$.  Suppose that 
\begin{enumerate}
  \item[i.]  $\tilde{L}^{L_k}_j \ge L_j$ and $\tilde{M}^{L_k}_\ell\le M_\ell$, for $j\ne k$ and $\ell=1,\ldots, m$, while
  \item[ii.]  $\tilde{L}^{M_k}_j\le L_j$ and $\tilde{M}^{M_k}_\ell \ge M_\ell$, $j=1\ldots, n$ and $\ell\ne k$.
\end{enumerate}
Then
\begin{equation}
d_{KS}(U,W)\leq ||\Delta f|| \left\{ \hbox{Var}(U) - (\lambda_1+\lambda_2)\right\} \enskip ,
\label{eq:bi.monotone.coupling}
\end{equation}
where $W\,\sim\, \hbox{Skellam}(\lambda_1,\lambda_2)$, with $\lambda_1, \lambda_2$ defined as in Theorem~\ref{thm:coupling}.
\label{thm:bi.monotone.coupling}
\end{theorem}
\noindent
The proof of this theorem is given in the appendix, in Section~\ref{sec:bi.monotone.coupling}.  We note that
for a collection of binary random variables to satisfy conditions $(i)$ and $(ii)$ in the above theorem, it is 
sufficient, for example, to generate a vector of positively associated random 
variables $(L_1,\ldots, L_n, 1-M_1,\ldots, 1-M_m)$.  The $L$'s and the $M$'s will then be positively associated
within, but negatively associated between, which in turn implies the conditions (i.e., analogous to 
positive and negative relatedness, respectively).  See~\citep[p. 78]{barbour2005introduction}, for example,
for a brief summary of these latter notions and their relationships.

With this theorem, the following then holds for large networks with dependent and homogeneous errors,
when the dependency is of the nature just defined.
\begin{corollary}
Suppose that the collections of edge indicator random variables $\{Y_{ij}\}_{\{i,j\}\in E^c}$ and
$\{Y_{ij}\}_{\{i,j\}\in E}$ satisfy conditions $(i)$ and $(ii)$ of Theorem~\ref{thm:bi.monotone.coupling} ,
playing the roles of the $L's$ and $M's$, respectively.  Then under assumptions (A1)-(A2),
\begin{equation}
d_{KS}\left(\, D_E\, ,\,  \hbox{Skellam}(\lambda,\lambda)\, \right) 
= O\left( \frac{\hbox{Var}(D_E) - 2\lambda}{2\lambda}\right) \enskip .
\label{eq:skellam.error.dep}
\end{equation}
\label{cor:skellam.error.dep}
\end{corollary}
This result can be compared to that of Theorem~\ref{thm:skellam.indep.case} , where the edge noise was independent and the error in approximating by a Skellam decayed like $n_v^{-1}$.  By way of contrast, Corollary~\ref{cor:skellam.error.dep} tells us that in order to achieve a decay in approximation error like $O(f(n_v))$, we must have $\hbox{Var}(D_E) = 2\lambda\left(1 + O(f(n_v))\right)$.

More generally, the quality of the approximation of $D_E$ by a Skellam will be influenced by the nature of the dependency in the errors, as the latter manifests itself through the overall variance $\hbox{Var}(D_E)$.  While the nature of that dependency is highly problem-specific, and a detailed investigation of possible cases is well beyond the scope of the present paper, nontrivial insight can be gained into the influence of the level of dependency on accuracy through numerical work under the following simple model.  

For a vector of binary random variables $(L_1,\ldots, L_n, 1-M_1,\ldots, 1-M_m)$, let $S=D + m$, where $D = \sum_{i=1}^n L_i - \sum_{i=1}^m M_i$.  We equip $S$ with a distribution of the form
\begin{equation}
\mathbb{P}\left(S = k\right) \propto {n+m \choose k}^{\nu - 1} 
                              \mathbb{P}\left(U + V = k\right) \enskip ,
\label{eq:COMB2}
\end{equation}
for $\nu$ a real number, where $U$ and $V$ are binomial, with parameters $(n,p)$ and $(m,q)$, respectively.  This distribution is a rescaling of that of the sum of two independent binomials, in a spirit analogous to the Conway-Maxwell binomial (COMB) distribution introduced recently by~\cite{kadane2014sums}.  The COMB distribution is a simple extension of the binomial distribution that introduces dependency among the corresponding Bernoulli random variables.  Our proposed distribution for $S$ in (\ref{eq:COMB2}) involves two binomial random variables, for which the corresponding Bernoulli random variables are dependent both within and between the two.  Accordingly, we call this a COMB2 distribution.

Now impose assumption (A2) on this model. Since the assumption implies that $\mathbb{E}\left[ D\right] = 0$, it follows that necessarily we must have $\mathbb{E}\left[ S\right] = m$.  Furthermore, the limiting Skellam distribution in Corollary~\ref{cor:skellam.error.dep} will be symmetric under this assumption.  Symmetry can be imposed here on the distribution of $S$, and hence $D$, by taking $n=m$ and $q=1-p$.  Therefore, we let $|E|={n_v\choose 2}/2 = n_v(n_v-1)/4$ and $\alpha = 1 - \beta$.  Note that this choice of $|E|$  means that our numerical work pertains to the case of dense graphs.  (We are unable to exhibit a sparse variant of the COMB2 model with the necessary characteristics above.)

Note that when $\nu=1$, the binary random variables $(L_1,\ldots, L_n, 1-M_1,\ldots, 1-M_m)$ are independent.  On the other hand, proceeding along lines of reasoning similar to those in~\cite{kadane2014sums}, it can be argued that the COMB2 distribution, with the parameter constraints just described, renders the $(L_1,\ldots, L_n, 1-M_1,\ldots, 1-M_m)$ positively associated when $\nu < 1$, with the mass being transferred increasingly to the endpoints of the support of the distribution of $D$ as $\nu\rightarrow - \infty$.  As a result, per the discussion immediately following Theorem~\ref{thm:bi.monotone.coupling}, the particular COMB2 distribution we have defined can be used to simulate network edge data in a way that satisfies the conditions of Corollary~\ref{cor:skellam.error.dep}.

\begin{figure}[bht]
\begin{center}
  \includegraphics[width=2.5in]{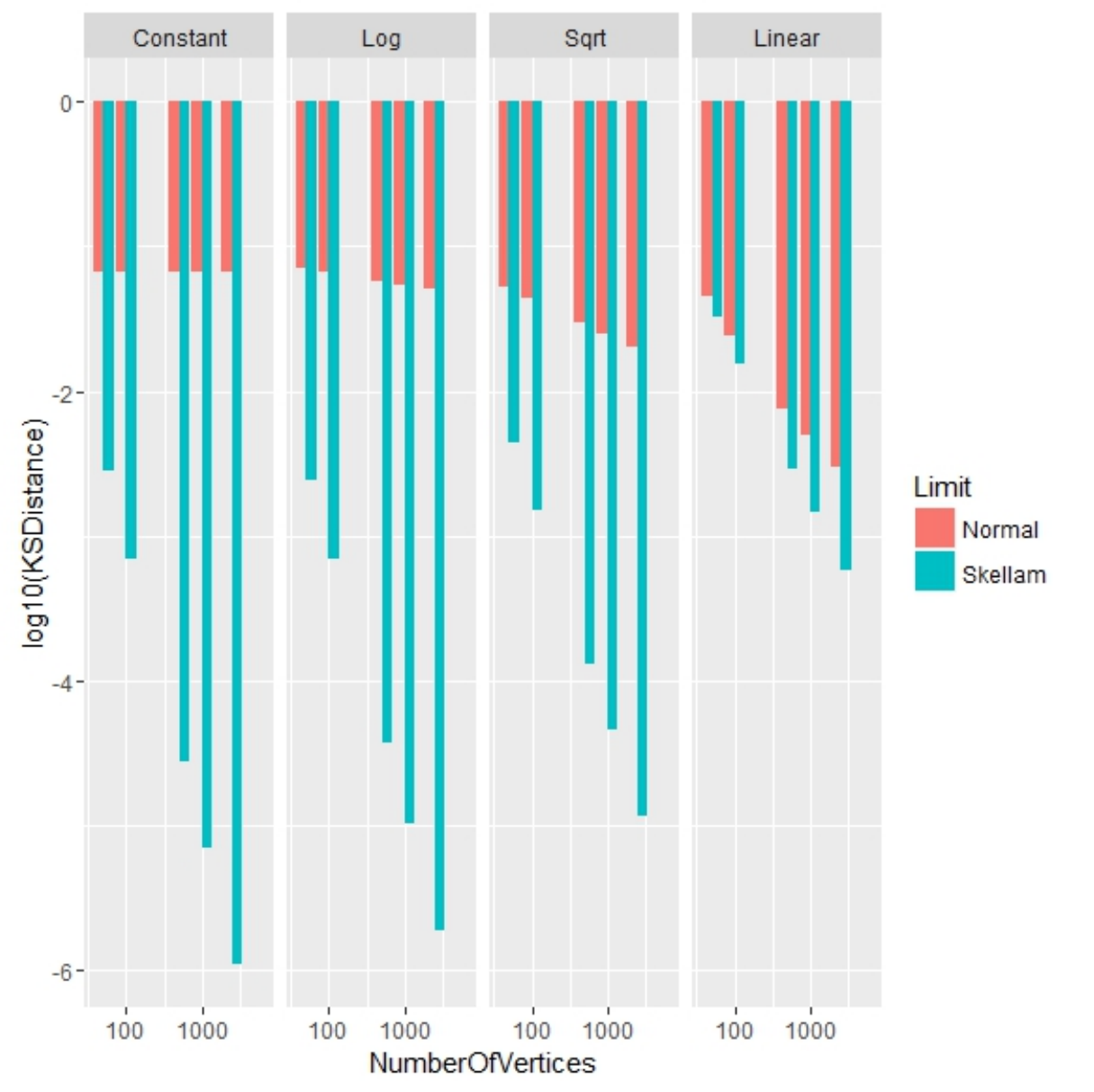}\quad
  \includegraphics[width=2.5in]{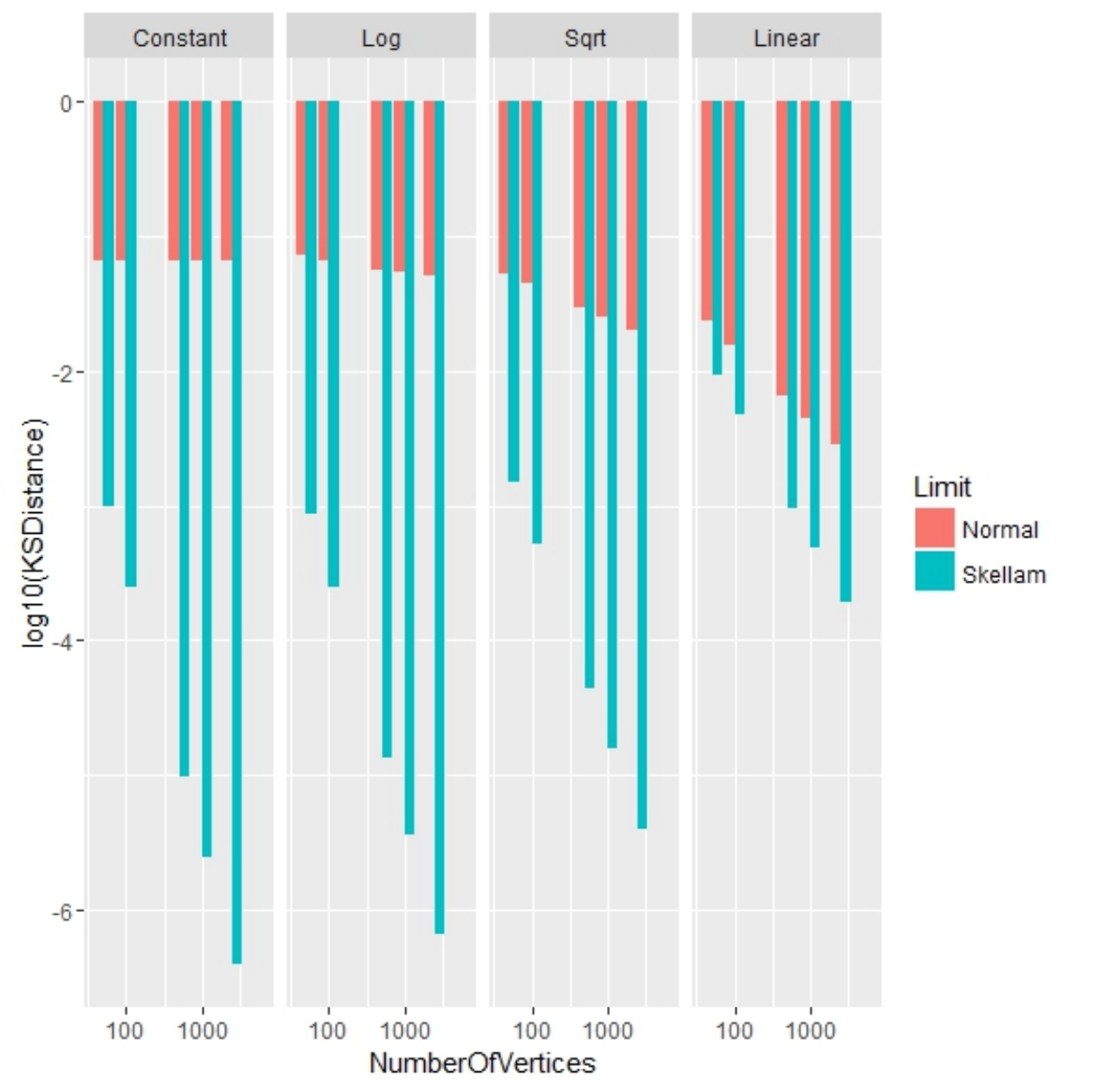}\\
  \includegraphics[width=2.5in]{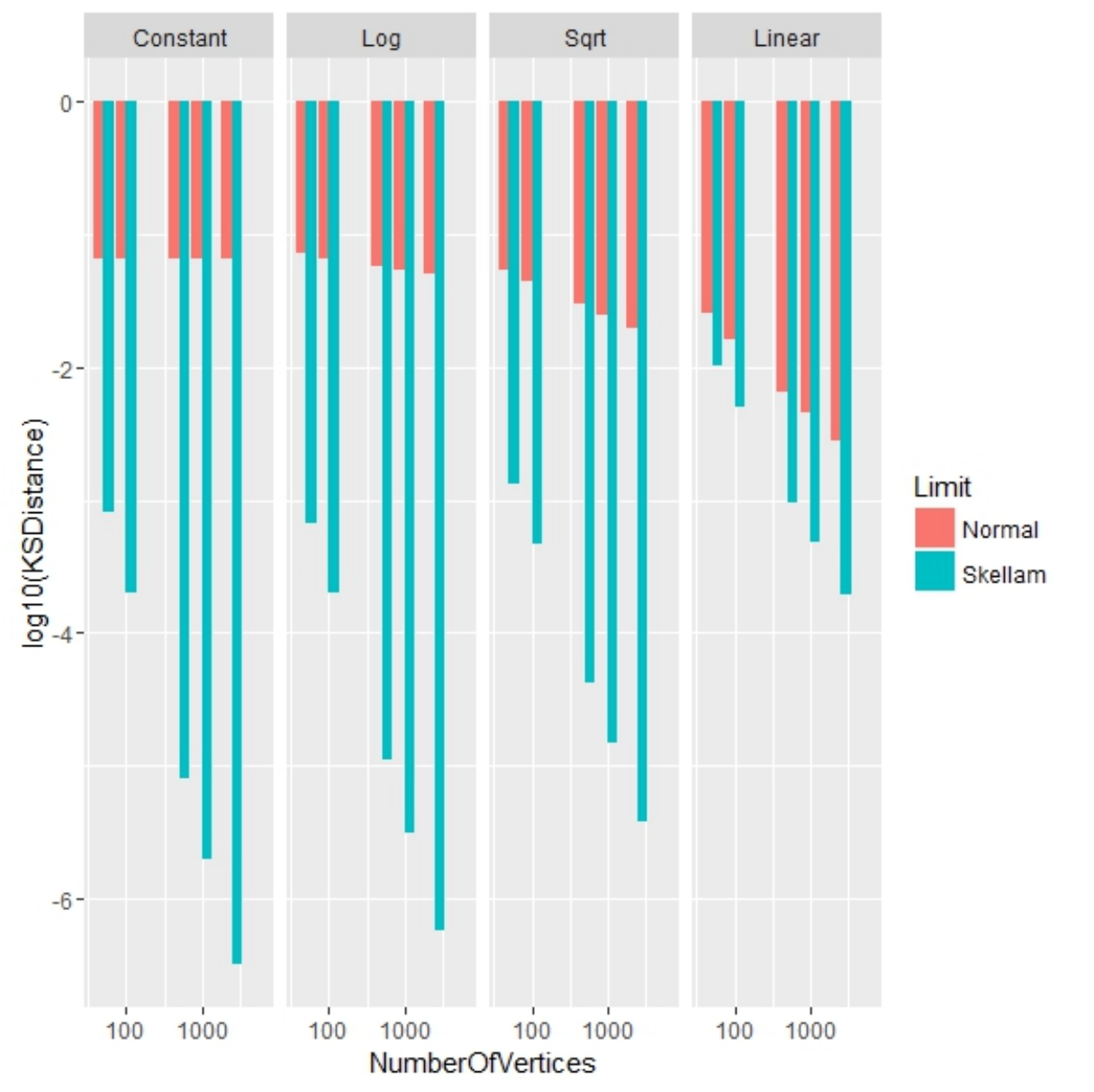}\quad
  \includegraphics[width=2.5in]{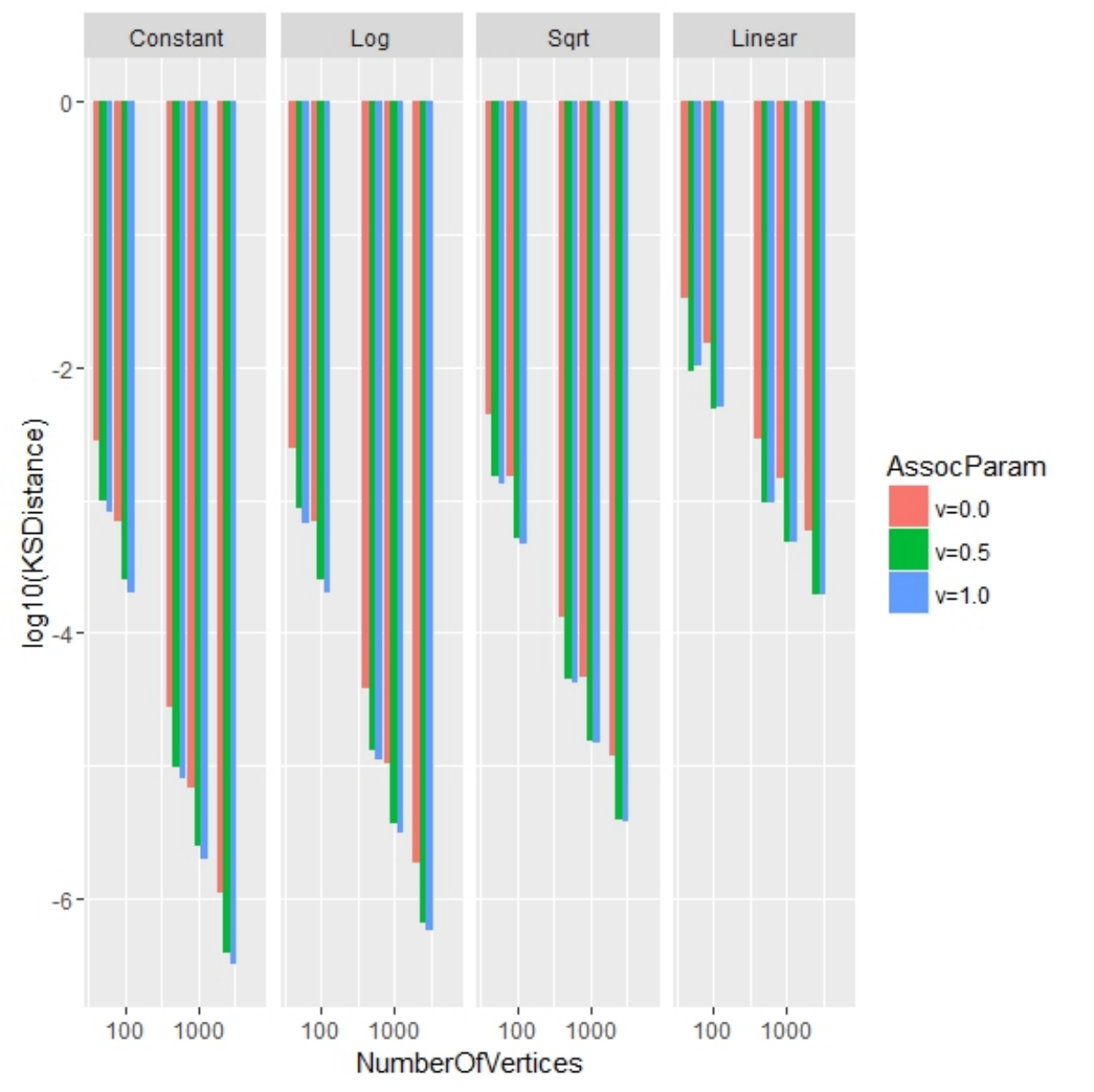}
\end{center}
\caption{(Log)Kolmogorov-Smirnov distances to the distribution of the discrepancy $D_E$ in edge counts under the COMB2 distribution, with $\nu=0$ (top left), $\nu=0.5$ (top right), and $\nu=1.0$ (bottom left), for the Skellam and standard normal approximations.  Also shown is a comparison of Skellam approximations as a function of $\nu$ (bottom right).}  
\label{fig:COMB2.simul}
\end{figure}

In Figure~\ref{fig:COMB2.simul} are shown the results of numerical work calculating the Kolmogorov-Smirnov distance between the Skellam and standard normal approximations to the distribution of the discrepancy $D_E$ in edge counts under the COMB2 distribution, for $\nu=0, 0.5$, and $1.0$.  The noise levels used here are the same as used earlier in producing Figure~\ref{fig:iid.simul}.  We see that the accuracy of the Skellam distribution decays slightly with increasing dependency in the errors, and with increasing noise levels.  

\section{Discussion}
\label{sec:discussion}

The propagation of uncertainty in network analysis is a topic that currently lags the field in development.  Despite almost 15 years of work in the modern `network science' era, on a vast array of topics, from researchers in many different disciplines, there remains a sizeable gap in our understanding of how `low level' errors (i.e., at the level of declaration of edge / non-edge status between vertex pairs) propagate to `high level' summaries (e.g., subgraph counts, 
centralities, etc.).  As a result, in most practical work, network summary statistics are cited without any indication
of likely error.    

Our contributions in this paper are aimed at helping to begin laying a foundation for work in this area, with a focus on  establishing an initial understanding of the distributional behavior of certain network summary statistics.  Our choice to work with subgraph counts is both natural and motivated by convenience, whereas our emphasis on the specific case of large networks with low-rate measurement error is intended to capture a sizeable fraction of what is encountered in practice.  Our formulation is reminiscent of the type of `signal plus noise' model commonly used in nonparametric function estimation and digital signal processing.  

In particular, in our formulation the true underlying graph $G$ is fixed.  This necessitates a different treatment than, say, traditional analysis of subgraph counts in classical Erdos-Renyi random graphs.  In the special case where an Erdos-Renyi model is assumed, as well as assuming independence among the measurement errors, and the analysis is done without conditioning on $G$, then the problem could be viewed as involving a classical random graph with two values for the probability of an edge arising in $\hat{G}$ (i.e., depending on whether or not there was an edge between a given pair of vertices in $G$).  In general, however, either when $G$ is fixed, as assumed in this paper, or from some other class of random graph models (e.g., various models with heterogeneous degree distributions), or when the measurement errors are dependent, the problem is more involved.  By conditioning on $G$, our formulation allows us to focus our analysis firstly on a high-level notion of Type I and II errors among (non)edges, and then secondly on the manner in which the structure of the underlying graph $G$ may interact with those errors.

We view our work as laying a key initial piece of the foundation on an important new problem area.  However, we have 
provided a detailed analysis only for the most fundamental of subgraph count statistics, i.e., the number of edges
in a network.  Our initial work on extension to counts for subgraphs of higher order suggests that the problem becomes increasingly nontrivial.  Specifically, the interaction of noise level, graph topology, and choice of subgraph would appear to need to be studied with care.  

The following general theorem should be useful in exploring further in this direction.  
\begin{theorem}
Let $H$ be a given subgraph of interest.  Re-express the difference $D$ in subgraph counts defined in equation (\ref{eq:subg.cnt.diff}) as
$$D_H = \sum_{{H'}\notin \mathcal{C}_H} L_{H'} 
                -  \sum_{{H'} \in \mathcal{C}_H} M_{H'} \enskip ,$$
for $\mathcal{C}_H = \{H' \subseteq K_{n_v}\, :\, H' \cong H , H' \subseteq G\}$,
where $L_{H'}$ and $M_{H'}$ are indicator variables of Type I and Type II error, respectively, for a subgraph $H'$.
Under the assumption of independent edge noise,
\begin{equation}
d_{KS}(D_H,W)\leq ||\Delta f|| \left\{ \hbox{Var}(D_H) - (\lambda_1+\lambda_2)\right\} \enskip ,
\label{eq:inter.bnd.gen.H}
\end{equation}
where $W\,\sim\, \hbox{Skellam}(\lambda_1,\lambda_2)$, with
$\lambda_1=\sum_{{H'}\notin \mathcal{C}_H} p_{H'}$ and 
$\lambda_2=\sum_{{H'} \in \mathcal{C}_H} q_{H'}$, for
$p_{H'} = \mathbb{E}\left[L_{H'}\right]$ and
$q_{H'} = \mathbb{E}\left[M_{H'}\right]$.
\label{thm:bnd.for.arb.H}
\end{theorem}
This result follows directly from application of Theorem~\ref{thm:bi.monotone.coupling} and the comment
immediately following that theorem.  In particular, each of the indicator random variables $L_{H'}$ and $1-M_{H'}$
may be expressed as a product of $n_{v(H)}$ choose two binary random variables, where $n_{v(H)}$
is the order of the subgraph $H$.  Since these products are non-decreasing functions of their arguments, 
and their arguments are independent, it follows that the collection of random variables defined by the union of the $L_{H'}$ and the $1-M_{H'}$ are positively associated (e.g., \citep{esary1967association}).

Application of this theorem to specific choices of subgraphs $H$ requires calculation or bounding of the two key quantities
within brackets in (\ref{eq:inter.bnd.gen.H}).  For the case of independent edge noise (which, nevertheless, yields {\em dependent} indicator variables $L_{H'}$ and $M_{H'}$), these quantities may be bounded through straightforward but tedious calculations for low-order subgraphs.  However, we also require control of the term $||\Delta f||$ in (\ref{eq:inter.bnd.gen.H}).  Under Conjecture~\ref{lipschitz:l1nel2}, this term is controlled by a term of order $(\lambda_1+\lambda_2)^{-1}$, but this conjecture, while supported by numerical work, remains to be proven.  

Our proof of Theorem~\ref{lipschitz}, in Section~\ref{sec:bounding.constant}, bounding $||\Delta f||$ for the case of $\lambda_1=\lambda_2$, required the control of alternating sums of differences of the ratios of modified Bessel functions of the first kind.  As such, the treatment is necessarily delicate.  Furthermore, the literature on quantities of this sort is lacking and, hence, we were required to develop several novel analysis results.  These results, which are of independent interest, are available in a separate manuscript~\citep{kash}.  Some extension thereof is presumably necessary to determine the validity of Conjecture~\ref{lipschitz:l1nel2}.

Finally, and interestingly, we mention that our preliminary numerical results suggest that in contexts like those of this paper, where the Skellam distribution obtains as the limiting distribution of the discrepancy in edge counts, an appropriate normal distribution may actually obtain for the discrepancy of counts of subgraphs of higher order, even for those of as little as order three (i.e., two-stars and triangles).  This observation suggests that what can be expected in this area are two regimes of limiting distributions, both normal and Poisson-like, in analogy to what is encountered in subgraph counting on classical random graphs, with the role of the Poisson distribution for results in this area replaced, in whole or in part, by the Skellam distribution.  

\section*{Acknowledgments}
We thank Kostas Spiliopoulos and Tasso Kaper for helpful discussions.
This work was supported in part by AFOSR award 12RSL042 and NSF grant CNS-0905565.

\section{Appendix}
\label{sec:appendix}

\subsection{Proof of Theorem~\ref{thm:steinop}.}
\label{sec:proof.of.steinop}

We begin with the operator,
\begin{equation} \nonumber
\mathcal{A}\left[f\left(k\right)\right] = \lambda_1f\left(k+1\right)-kf\left(k\right)-\lambda_2f\left(k-1\right)
\end{equation}
with the intent of showing that the random variable $W\sim\mbox{Skellam}\left(\lambda_1,\lambda_2\right)$ if and only if $\mathbb{E}\mathcal{A}\left[f\left(W\right)\right]=0$ for any bounded function $f:\mathbb{Z}\mapsto\mathbb{R}$.

We begin with the necessity direction and the computation of
\begin{eqnarray} \nonumber
\mathbb{E}\mathcal{A}\left[f\left(W\right)\right] &=& \mathbb{E}\left[\lambda_1f\left(W+1\right)-Wf\left(W\right)-\lambda_2f\left(W-1\right)\right]\\
 &\propto& \nonumber \sum_{k=-\infty}^\infty\left[\lambda_1f\left(k+1\right)-kf\left(k\right)-\lambda_2f\left(k-1\right)\right]\left(\sqrt{\frac{\lambda_1}{\lambda_2}}\right)^kI_k
\end{eqnarray}
where $\propto$ is to be read as ``proportional to,'' and as shorthand, we write $I_k$ for $I_k(2\sqrt{\lambda_1\lambda_2})$.  
By standard properties of $I_k$ (e.g., \citep{abramowitz1972handbook}) we have that
\begin{equation} \nonumber
I_{k-1}-I_{k+1}=\frac{k}{\sqrt{\lambda_1\lambda_2}}I_k
\end{equation}
or, in other words,
\begin{equation}\label{recur}
\sqrt{\lambda_1\lambda_2}\frac{I_{k-1}}{I_k}-\sqrt{\lambda_1\lambda_2}\frac{I_{k+1}}{I_k}=k.
\end{equation}
This means that
\begin{eqnarray} \nonumber
\mathbb{E}\mathcal{A}\left[f\left(W\right)1\left\{W\leq n\right\}\right] &\propto& \nonumber \sum_{k=-\infty}^n\left[\sqrt{\frac{\lambda_1}{\lambda_2}}f\left(k+1\right)-\frac{I_{k-1}}{I_k}f\left(k\right)+\right.\\
 & & \nonumber \left.\frac{I_{k+1}}{I_k}f\left(k\right)-\sqrt{\frac{\lambda_2}{\lambda_1}}f\left(k-1\right)\right]\left(\sqrt{\frac{\lambda_1}{\lambda_2}}\right)^kI_k\\
 &=& \nonumber \sum_{k=-\infty}^n\left[\left(\sqrt{\frac{\lambda_1}{\lambda_2}}\right)^{k+1}I_kf\left(k+1\right)-\left(\sqrt{\frac{\lambda_1}{\lambda_2}}\right)^kI_{k-1}f\left(k\right)\right]\\
  & & \nonumber +\sum_{k=-\infty}^n\left[\left(\sqrt{\frac{\lambda_1}{\lambda_2}}\right)^kI_{k+1}f\left(k\right)-\left(\sqrt{\frac{\lambda_1}{\lambda_2}}\right)^{k-1}I_kf\left(k-1\right)\right]\\
 &=& \nonumber \left(\sqrt{\frac{\lambda_1}{\lambda_2}}\right)^{n+1}I_nf\left(n+1\right)+\left(\sqrt{\frac{\lambda_1}{\lambda_2}}\right)^nI_{n+1}f\left(n\right).
\end{eqnarray}
Now, since $f$ is bounded,
\begin{equation} \nonumber
\lim_{n\rightarrow\infty}\left\{\left(\sqrt{\frac{\lambda_1}{\lambda_2}}\right)^{n+1}I_nf\left(n+1\right)+\left(\sqrt{\frac{\lambda_1}{\lambda_2}}\right)^nI_{n+1}f\left(n\right)\right\}=0
\end{equation}
so that by monotone convergence,

$$
\begin{aligned}
\mathbb{E}\left[\mathcal{A}\left[f(W)\right] \right] &= \lim_{n\rightarrow \infty} \mathbb{E}\mathcal{A}\left[f\left(W\right)1\left\{W\leq n\right\}\right]\\
& =\lim_{n\rightarrow\infty}\left\{\left(\sqrt{\frac{\lambda_1}{\lambda_2}}\right)^{n+1}I_nf\left(n+1\right)+\left(\sqrt{\frac{\lambda_1}{\lambda_2}}\right)^nI_{n+1}f\left(n\right)\right\}\\
&=0
\end{aligned}
$$
which proves the claim.

To prove sufficiency, we begin with $\mathbb{E}\mathcal{A}\left[f\left(W\right)\right]=0$ and suppose that $f_k\left(j\right)=1\left\{j=k\right\}$ for some $j\in\mathbb{Z}$ in which case
\begin{equation} \nonumber \lambda_1p\left(k-1\right)-kp\left(k\right)-\lambda_2p\left(k+1\right)=0
\end{equation}
where $p\left(k\right)=\mathbb{P}\left(W=k\right)$.  An ansatz of
\begin{equation} \nonumber
\begin{array}{ccc}
S\left(k\right) = \left(\sqrt{\frac{\lambda_1}{\lambda_2}}\right)^kI_k\left(2\sqrt{\lambda_1\lambda_2}\right) & \mbox{and} & T\left(k\right) = \left(\sqrt{\frac{\lambda_1}{\lambda_2}}\right)^kK_k\left(2\sqrt{\lambda_1\lambda_2}\right)
\end{array}
\end{equation}
shows that $S$ and $T$ form two linearly independent solutions to this second order linear difference equation, where  $I_k(x)$ and $K_k(x)$ are the modified Bessel functions of the first and second kinds.  Thus, we know that the general solution is given by,
\begin{equation} \nonumber
p\left(k\right) = C_1S\left(k\right)+C_2T\left(k\right)
\end{equation}
for some constants $C_1,C_2\in\mathbb{R}$.   

Now, to determine the constants $C_1$ and $C_2$ we appeal to the fact that $\sum_{k=-\infty}^\infty p\left(k\right)=1$.  Since $I_k,K_k>0$ for all $k\in\mathbb{Z}$ and $\sum_{k=-\infty}^\infty K_k=\infty$ it must be that $C_2=0$.  Now, consider the generating function
\begin{equation} \nonumber
e^{\frac{z}{2}\left(t+1/t\right)}=\sum_{k=-\infty}^\infty t^kI_k\left(z\right)
\end{equation}
which means that
\begin{eqnarray} \nonumber
C_1 &=& \frac{1}{\sum_{k=-\infty}^\infty\left(\sqrt{\frac{\lambda_1}{\lambda_2}}\right)^kI_k\left(2\sqrt{\lambda_1\lambda_2}\right)}\\
 &=& \nonumber \frac{1}{e^{\sqrt{\lambda_1\lambda_2}\left(\sqrt{\frac{\lambda_1}{\lambda_2}}+\sqrt{\frac{\lambda_2}{\lambda_1}}\right)}}\\
 &=&  \nonumber e^{-\left(\lambda_1+\lambda_2\right)}
\end{eqnarray}
so that
\begin{equation} \nonumber
p\left(k\right) = e^{-\left(\lambda_1+\lambda_2\right)}\left(\sqrt{\frac{\lambda_1}{\lambda_2}}\right)^kI_k\left(2\sqrt{\lambda_1\lambda_2}\right)
\end{equation}
so that $W\sim\mbox{Skellam}\left(\lambda_1,\lambda_2\right)$.

\begin{flushright}
$\square$
\end{flushright}

\subsection{Proof of Theorem~\ref{thm:coupling}.}
\label{sec:proof.of.coupling}

Given that $f_x$ is a solution to $\mathcal{A}[f_x(k)]=g_x(k)$, we have
$$\lambda_1f_x(k+1)-kf_x(k)-\lambda_2 f_x(k-1) = 1(k\leq x) - \mathbb{P}\left[W\leq x\right].$$

Substituting $k=U$ and taking expected values, we obtain,
\begin{equation}
\label{A3.1}
\begin{aligned}
\left| \mathbb{P}\left[ U\leq x\right] - \mathbb{P}\left[ W\leq x\right] \right| = \left| \mathbb{E}\left[ \lambda_1 f_x(U+1)-Uf_x(U)-\lambda_2 f_x(U-1) \right]\right|\enskip.
\end{aligned}
\end{equation}

Next, recall from (\ref{V Def}) that $U=\sum_{k=1}^n L_k -\sum_{k=1}^m M_k$.  Since $\lambda_1=\sum_{k=1}^n p_k$ and $\lambda_2=\sum_{k=1}^m q_k$, we have after conditioning on $L_k$ and $M_k$,

$$
\begin{aligned}
& \left| \mathbb{E}\left[ \lambda_1 f_x(U+1)-Uf_x(U)-\lambda_2 f_x(U-1)\right] \right|\\
&= \left| \sum_{k=1}^n \mathbb{E}\left[p_k f_x(U+1) - L_kf_x(U) \right] + \sum_{k=1}^m \mathbb{E}\left[M_kf_x(U)-q_kf_x(U-1)\right] \right|\\
&=\left|\sum_{k=1}^n p_k \left(\mathbb{E}\left[ f_x(U+1)\right]-\mathbb{E}\left[f_x(U)|L_k=1\right] \right)+ \sum_{k=1}^m q_k\left( \mathbb{E}\left[f_x(U)|M_k=1\right]-\mathbb{E}\left[f_x(U-1)\right]\right)\right|\\
&=\left|\sum_{k=1}^n p_k \left( \mathbb{E}\left[ f_x(U+1)-f_x\left(U_k^{(L)}+1\right) \right]\right)+\sum_{k=1}^m q_k \left( \mathbb{E}\left[ f_x\left(U_k^{(M)}-1\right) - f_x(U-1)\right]\right)\right|\\
&\leq \sum_{k=1}^n p_k ||\Delta f|| \mathbb{E}\left|U-U_k^{(L)}\right| + \sum_{k=1}^m q_k ||\Delta f|| \mathbb{E}\left| U-U_k^{(M)}\right|\\
& = ||\Delta f|| \left[\sum_{k=1}^n p_k   \mathbb{E}\left|U-U_k^{(L)}\right| + \sum_{k=1}^m q_k  \mathbb{E}\left| U-U_k^{(M)}\right| \right].\\
\end{aligned}
$$
Combining this with (\ref{A3.1}) yields the result.

\begin{flushright}
$\square$
\end{flushright}

\subsection{Proof of Theorem~\ref{steinsol}.}
\label{sec:proof.of.steinsol}

First, consider the solution to  
\begin{equation} 
\label{Stein Equation Solution}
\lambda_1f\left(k+1\right)-kf\left(k\right)-\lambda_2f\left(k-1\right)=g\left(k\right),
\end{equation}
for some bounded function $g:\mathbb{Z}\mapsto\mathbb{R}$, with the boundary condition
\begin{equation}\label{boundary}
\lim_{k\rightarrow-\infty}\left(\sqrt{\frac{\lambda_1}{\lambda_2}}\right)^kI_kf\left(k\right)=0.
\end{equation}
We use (\ref{recur}) to substitute for $k$ in (\ref{Stein Equation Solution}).  
Then, multiplying both sides of (\ref{Stein Equation Solution}) by $\left(\sqrt{{\lambda_1}/{\lambda_2}}\right)^kI_k$, we obtain,  
 $$
 \begin{aligned}
& \lambda_1 \left( \sqrt{\frac{\lambda_1}{\lambda_2}} \right)^k I_kf(k+1)-\lambda_1\left( \sqrt{\frac{\lambda_1}{\lambda_2}} \right)^{k-1} I_{k-1}f(k)\\
&+\lambda_2\left( \sqrt{\frac{\lambda_1}{\lambda_2}} \right)^{k+1} I_{k+1}f(k)-\lambda_2\left( \sqrt{\frac{\lambda_1}{\lambda_2}} \right)^k I_kf(k-1)\\
&=\left(\sqrt{\frac{\lambda_1}{\lambda_2}}\right)^kI_k g(k)
\end{aligned}
 $$
 which is the same as,
 
  $$
 \begin{aligned}
& \left( \sqrt{\frac{\lambda_1}{\lambda_2}} \right)^{k+1} I_kf(k+1)-\left( \sqrt{\frac{\lambda_1}{\lambda_2}} \right)^{k} I_{k-1}f(k)\\
&+\left( \sqrt{\frac{\lambda_1}{\lambda_2}} \right)^{k} I_{k+1}f(k)-\left( \sqrt{\frac{\lambda_1}{\lambda_2}} \right)^{k-1} I_kf(k-1)\\
&=\frac{1}{\sqrt{\lambda_1\lambda_2}}\left(\sqrt{\frac{\lambda_1}{\lambda_2}}\right)^kI_k g(k).
\end{aligned}
 $$
 
 Notice that we have grouped terms together so that summing over $k$ yields a telescoping sum.  So, summing over $k\in\left\{-\infty,\ldots,n\right\}$ and using the boundary condition (\ref{boundary}), 
\begin{equation} \nonumber
\left(\sqrt{\frac{\lambda_1}{\lambda_2}}\right)^{n+1}I_nf\left(n+1\right)+\left(\sqrt{\frac{\lambda_1}{\lambda_2}}\right)^nI_{n+1}f\left(n\right)=\frac{1}{\sqrt{\lambda_1\lambda_2}}\sum_{k=-\infty}^n\left(\sqrt{\frac{\lambda_1}{\lambda_2}}\right)^kI_kg\left(k\right).
\end{equation}
Now, multiplying both sides by $\left(-1\right)^{n+1}/(I_nI_{n+1})$ and summing over $n\in\left\{c,c+1,\ldots,m\right\}$ for $m>c$ and over $n\in\left\{m,m+1,\ldots,c-1\right\}$, for some initial condition $c\in\mathbb{Z}$ and $f\left(c\right)\in\mathbb{R}$, we obtain
\begin{equation} \nonumber
f\left(m\right) = \left\{\begin{array}{ll}
\left(-1\right)^m\left(\sqrt{\frac{\lambda_2}{\lambda_1}}\right)^mI_m\left[\left(-1\right)^c\left(\sqrt{\frac{\lambda_1}{\lambda_2}}\right)^c\frac{1}{I_c}f\left(c\right)\right. & \\
\left.+\frac{1}{\sqrt{\lambda_1\lambda_2}}\sum_{n=c}^{m-1}\frac{\left(-1\right)^{n+1}}{I_nI_{n+1}}\sum_{k=-\infty}^n\left(\sqrt{\frac{\lambda_1}{\lambda_2}}\right)^kI_kg\left(k\right)\right] & \mbox{if $m>c$}\\
\left(-1\right)^m\left(\sqrt{\frac{\lambda_2}{\lambda_1}}\right)^mI_m\left[\left(-1\right)^c\left(\sqrt{\frac{\lambda_1}{\lambda_2}}\right)^c\frac{1}{I_c}f\left(c\right)\right. & \\
\left.-\frac{1}{\sqrt{\lambda_1\lambda_2}}\sum_{n=m}^{c-1}\frac{\left(-1\right)^{n+1}}{I_nI_{n+1}}\sum_{k=-\infty}^n\left(\sqrt{\frac{\lambda_1}{\lambda_2}}\right)^kI_kg\left(k\right)\right] & \mbox{if $m<c$}.
\end{array}\right.
\end{equation}

Note that if
\begin{equation} \nonumber
g\left(k\right) = g_x(k)= 1\left\{k\leq x\right\}-\mathbb{P}\left(W\leq x\right)
\end{equation}
then
\begin{equation} \nonumber
\sum_{k=-\infty}^n\left(\sqrt{\frac{\lambda_1}{\lambda_2}}\right)^kI_kg\left(k\right) = \left\{\begin{array}{ll}
e^{\lambda_1+\lambda_2}\mathbb{P}\left(W\leq n\right)\mathbb{P}\left(W>x\right) & \mbox{if $n\leq x$}\\
e^{\lambda_1+\lambda_2}\mathbb{P}\left(W\leq x\right)\mathbb{P}\left(W>n\right) & \mbox{if $n\geq x$}
\end{array}\right.
\end{equation}
since, for example if $n\leq x$
\begin{eqnarray} \nonumber
\sum_{k=-\infty}^n\left(\sqrt{\frac{\lambda_1}{\lambda_2}}\right)^kI_kg\left(k\right) &=& e^{\lambda_1+\lambda_2}\sum_{k=-\infty}^n\mathbb{P}\left(W=k\right)g\left(k\right)\\
 &=& \nonumber e^{\lambda_1+\lambda_2}\sum_{k=-\infty}^n\mathbb{P}\left(W=k\right)\left[1\left\{k\leq x\right\}-\mathbb{P}\left(W\leq x\right)\right]\\
 &=& \nonumber e^{\lambda_1+\lambda_2}\left[\mathbb{P}\left(W\leq\min\left\{n,x\right\}\right)-\mathbb{P}\left(W\leq x\right)\mathbb{P}\left(W\leq n\right)\right]\\
 &=&  \nonumber e^{\lambda_1+\lambda_2}\mathbb{P}\left(W\leq n\right)\left[1-\mathbb{P}\left(W\leq x\right)\right]\\
 &=& \nonumber e^{\lambda_1+\lambda_2}\mathbb{P}\left(W\leq n\right)\mathbb{P}\left(W>x\right).
\end{eqnarray}
The case that $n\geq x$ is similar.  This means that
\begin{equation} \nonumber
f_x\left(m\right) = \left\{\begin{array}{ll}
\left(-1\right)^m\left(\sqrt{\frac{\lambda_2}{\lambda_1}}\right)^mI_m\left[\left(-1\right)^c\left(\sqrt{\frac{\lambda_1}{\lambda_2}}\right)^c\frac{1}{I_c}f\left(c\right)\right. & \\
\left.+\frac{e^{\lambda_1+\lambda_2}}{\sqrt{\lambda_1\lambda_2}}\sum_{n=c}^{m-1}\frac{\left(-1\right)^{n+1}}{I_nI_{n+1}}\mathbb{P}\left(W\leq\min\left\{n,x\right\}\right)\mathbb{P}\left(W>\max\left\{n,x\right\}\right)\right] & \mbox{if $m>c$}\\
\left(-1\right)^m\left(\sqrt{\frac{\lambda_2}{\lambda_1}}\right)^mI_m\left[\left(-1\right)^c\left(\sqrt{\frac{\lambda_1}{\lambda_2}}\right)^c\frac{1}{I_c}f\left(c\right)\right. & \\
\left.-\frac{e^{\lambda_1+\lambda_2}}{\sqrt{\lambda_1\lambda_2}}\sum_{n=m}^{c-1}\frac{\left(-1\right)^{n+1}}{I_nI_{n+1}}\mathbb{P}\left(W\leq\min\left\{n,x\right\}\right)\mathbb{P}\left(W>\max\left\{n,x\right\}\right)\right] & \mbox{if $m<c$}.
\end{array}\right.
\end{equation}

\begin{flushright}
$\square$
\end{flushright}

\subsection{Proof of Theorem~\ref{lipschitz}.}
\label{sec:proof.of.lipschitz}

Our proof of Theorem~\ref{lipschitz} is highly involved, from an analysis perspective, but the overall program can be stated
in a relatively succinct manner.  Accordingly, we sketch here the overall program behind our proof and refer the interested reader
to the Supplementary Materials for a detailed account.
 
Recall that we are trying to obtain a bound on $|\Delta f_x(j)|=|f_x(j+1)-f_x(j)|$ independent of $x\in \mathbb{R}$ and $j\in \mathbb{Z}$.  From Theorem~\ref{steinsol}, we have the solution to the Stein equation, however to use it to bound $|\Delta f_x(j)|$, we need to simplify it further.  For ease of notation, we simply refer to $f$ instead of $f_x$ and $g$ instead of $g_x$.

First, note that we have the freedom to choose the initial condition $(c,f(c))$.  Making the choice
that $c=\lambda_2-\lambda_1$, and hence that $c=0$ under the assumption that $\lambda_1=\lambda_2$,
we are able to simplify our expression for $f$ in Theorem~\ref{steinsol} to read, in the case that  $m>0$, as
\begin{equation}
\label{split1}
\begin{split}
f(m+1)  &= \left(-1\right)^{m+1}I_{m+1}\Bigg[\frac{1}{I_0}f\left(0\right)\\
 & \qquad -\frac{e^{2\lambda}}{\lambda}\frac{1}{I_0I_1}\mathbb{P}\left(W\leq\min\left\{0,x\right\}\right)\mathbb{P}\left(W>\max\left\{0,x\right\}\right)\\
 & \qquad
 -\frac{e^{2\lambda}}{\lambda}\sum_{n=0}^{m-1}\frac{\left(-1\right)^{n+1}}{I_{n+1}I_{n+2}}\mathbb{P}\left(W\leq\min\left\{n+1,x\right\}\right)\mathbb{P}\left(W>\max\left\{n+1,x\right\}\right)\Bigg]
\end{split}
\end{equation}
and, in the case that if $m<0$, as
\begin{equation*}
\begin{split}
f(m-1)  &= \left(-1\right)^{m-1}I_{m-1}\Bigg[\frac{1}{I_0}f\left(0\right)\\
 & \qquad -\frac{e^{2\lambda}}{\lambda}\frac{1}{I_0I_1}\mathbb{P}\left(W\leq\min\left\{-1,x\right\}\right)\mathbb{P}\left(W>\max\left\{-1,x\right\}\right)\\
 & \qquad
 +\frac{e^{2\lambda}}{\lambda}\sum_{n=m}^{-1}\frac{\left(-1\right)^{n+1}}{I_{n-1}I_n}\mathbb{P}\left(W\leq\min\left\{n-1,x\right\}\right)\mathbb{P}\left(W>\max\left\{n-1,x\right\}\right)\Bigg].
\end{split}
\end{equation*}
Finally, for the case $m=0$, we have
\begin{eqnarray}
f\left(0\right)  &=& \label{f0 Definition}  \frac{e^{2\lambda}}{2\lambda}\frac{1}{I_0+I_1}\left[\mathbb{P}\left(W\leq\min\left\{0,x\right\}\right)\mathbb{P}\left(W>\max\left\{0,x\right\}\right)\right.\\
  & & \nonumber \left.+\mathbb{P}\left(W\leq\min\left\{-1,x\right\}\right)\mathbb{P}\left(W>\max\left\{-1,x\right\}\right)\right].
\end{eqnarray}

Next, through manipulation of the arguments in the sums defining the above expressions for $f$, exploiting properties of the
modified Bessel functions $I_k$, and applying the triangle inequality, we are able to produce bounds on the differences
$|f(m+1) - f(m)|$ of the form
\begin{eqnarray}\label{majorcase}
\\ \nonumber
\left|f\left(m+1\right)-f\left(m\right)\right| &\leq & \frac{\mathbb{P}\left(W\leq x\right)}{\lambda}\left\{\left|\sum_{n=1,3,\ldots}^{m-1}\frac{I_{m+1}}{I_{n+2}}-\frac{I_m}{I_{n-1}}\right|\right.\\
 & & \nonumber \left. +\sum_{n=1,3,\ldots}^{m-1}H\left(n\right)\left|\frac{I_{m+1}I_n}{I_{n+1}I_{n+2}}-\frac{I_m}{I_{n+1}}-\frac{I_{m+1}}{I_{n+1}}+\frac{I_m}{I_{n-1}}\right| \right\}\\
 & & \nonumber +\left|I_m\frac{1}{I_0}f\left(0\right)+I_{m+1}\left\{\frac{1}{I_0}f\left(0\right)\right.\right. \\
  & & \nonumber \left.\left.-\frac{e^{2\lambda}}{\lambda}\frac{1}{I_0I_1}\mathbb{P}\left(W\leq\min\left\{0,x\right\}\right)\mathbb{P}\left(W>\max\left\{0,x\right\}\right)\right\}\right| \enskip ,
\end{eqnarray}
 if $m$ is even, and
\begin{eqnarray}\label{majorcase2}\\ \nonumber
\left|f\left(m+1\right)-f\left(m\right)\right| &\leq & \frac{\mathbb{P}\left(W\leq x\right)}{\lambda}\left\{\left|\sum_{n=1,3,\ldots}^{m-1}\frac{I_{m+1}}{I_{n+2}}-\frac{I_m}{I_{n-1}}\right|\right.\\
 & & \nonumber +\sum_{n=1,3,\ldots}^{m-1}H\left(n\right)\left|\frac{I_{m+1}I_n}{I_{n+1}I_{n+2}}-\frac{I_m}{I_{n+1}}-\frac{I_{m+1}}{I_{n+1}}+\frac{I_m}{I_{n-1}}\right| \\
  & & \nonumber \left. +\left|H\left(m+1\right)-H\left(m\right)\right|\right\}\\
 & & \nonumber +\left|I_m\frac{1}{I_0}f\left(0\right)+I_{m+1}\left\{\frac{1}{I_0}f\left(0\right)\right.\right.\\
  & & \nonumber \left.\left.-\frac{e^{2\lambda}}{\lambda}\frac{1}{I_0I_1}\mathbb{P}\left(W\leq\min\left\{0,x\right\}\right)\mathbb{P}\left(W>\max\left\{0,x\right\}\right)\right\}\right| \enskip ,
\end{eqnarray}
if $m$ is odd.  Here $H(n) = \mathbb{P}\left(W > n\right)/\mathbb{P}\left(W=n\right)$ is the inverse of the
hazard function of the Skellam distribution (and is not to be confused with our use of $H$ in the main body of the paper
as a subgraph of the graph $G$).

Note that (\ref{majorcase}) is defined by three key terms, while (\ref{majorcase2})
has the same three,  augmented by the addition of a fourth, i.e., $\left|H\left(m+1\right)-H\left(m\right)\right|$.  
Through a series of arguments (the result for  each of which is presented as a separate proposition in the Supplementary Materials),
we are able to control each of these terms as follows.  First, we show that
\begin{equation}
\sup_{m\in\mathbb{N}^+}\left|\sum_{n=1,3,\ldots}^{m-1}\frac{I_{m+1}}{I_{n+2}}-\frac{I_m}{I_{n-1}}\right|\leq 5\enskip .
\end{equation}
Next we show that
\begin{equation}
 \sum_{n=1,3,\ldots}^{m-1}H\left(n\right)\left|\frac{I_{m+1}I_n}{I_{n+1}I_{n+2}}-\frac{I_m}{I_{n+1}}-\frac{I_{m+1}}{I_{n+1}}+\frac{I_m}{I_{n-1}}\right| \leq 73 \enskip ,
\end{equation}
for $\lambda\geq1$.
And furthermore, we show that
\begin{eqnarray}
& & \nonumber \left|I_m\frac{1}{I_0}f\left(0\right)+I_{m+1}\left\{\frac{1}{I_0}f\left(0\right)-\frac{e^{2\lambda}}{\lambda}\frac{1}{I_0I_1}\mathbb{P}\left(W\leq\min\left\{0,x\right\}\right)\mathbb{P}\left(W>\max\left\{0,x\right\}\right)\right\}\right|\\
 &\le& \nonumber \frac{\mathbb{P}\left(W\leq  x \right)}{\lambda} \enskip. 
\end{eqnarray}
Finally, it is clear that
\begin{eqnarray}
& & \nonumber H\left(m\right)-H\left(m+1\right)\\
 &=& \nonumber \frac{\mathbb{P}\left(W>m\right)}{\mathbb{P}\left(W=m\right)}-\frac{\mathbb{P}\left(W>m+1\right)}{\mathbb{P}\left(W=m+1\right)}\\
 &=& \nonumber \frac{1}{\mathbb{P}\left(W=m\right)}\left[\mathbb{P}\left(W>m\right)-\mathbb{P}\left(W>m+1\right)\right]\\
 &=& \nonumber \frac{\mathbb{P}\left(W=m+1\right)}{\mathbb{P}\left(W=m\right)}\\
 &\leq& \nonumber 1
\end{eqnarray}
and, at the same time $H\left(m\right)-H\left(m+1\right)\geq 0$ so, we have that we may bound the magnitude of 
this final term by $1$.

As a result of all of the above, we may conclude that
\begin{equation}
\left|f\left(m+1\right)-f\left(m\right)\right|\le  \frac{80}{\lambda}
\end{equation}
for $m>0$.  Or, equivalently, we may express the right-hand side above as $160/2\lambda$.

The argument for the case of $m<0$ involves similar reasoning, as described in the Supplementary Materials.
\begin{flushright}
$\square$
\end{flushright}

\subsection{Proof of Theorems~\ref{thm:skellam.indep.case} and~\ref{thm:normal.indep.case}.}
\label{sec:proof.of.normal.vs.skellam.indep}

\subsubsection{Proof of Theorem~\ref{thm:skellam.indep.case}.}

The terms $\mathbb{E}\left|U-U_k^{(L)} \right|$ and $\mathbb{E}\left| U-U_k^{(M)}\right|$ 
in (\ref{Stein KS Bound}) measure the dependence of $U$ on the events $L_k=1$ and $M_k=1$, respectively.  
In the context of the empirical graph $\hat{G}$, the random variables $L$ are equal to $Y_{ij}$, for 
$\{i,j\}\in E^c$, while the random variables $M$ are equal to $Y_{ij}$, for $\{i,j\}\in E$.
With the $Y_{ij}$ assumed independent, $U_k^{(L)}$ and $U_k^{(M)}$ 
are independent of their respective events, and so we obtain
\begin{equation}
d_{KS}(U,W)\leq ||\Delta f|| \left [\sum_{k=1}^n p_k^2 + \sum_{k=1}^m q_k^2\right] \enskip .
\label{eq:stein.bnd.indep}
\end{equation}

Accordingly, and drawing on definitions and the result of Theorem~\ref{lipschitz},
$$
\begin{aligned}
d_{KS}\left(D_E, \hbox{Skellam}(\lambda,\lambda)\right)\leq \frac{1}{\lambda} \left[\sum_{(i,j)\in E^c} \alpha^2+\sum_{(i,j)\in E} \beta^2\right]&=\frac{|E^c| \alpha^2 + |E|\beta^2}{|E^c|\alpha}\\
&=\frac{|E^c| \alpha^2 + |E|\left(\frac{|E^c|}{|E|}\right)^2\alpha^2}{|E^c|\alpha}\\
&=\alpha +\frac{|E^c|}{|E|}\alpha\\
&=\alpha + \frac{{n_v\choose 2}-|E|}{|E|}\alpha\\
&=\frac{{n_v\choose 2}}{|E|}\alpha \enskip .
\end{aligned}
$$
Noting that $\alpha = \lambda / |E^c|$, and recalling that $|E^c|=\Theta\left(n_v^2\right)$ under both sparse and dense graphs $G$, the last quantity above is seen to behave like $\lam/|E|$ which, under
assumption (A3) and our definition of sparse and dense in Section~\ref{sec:initial.results} , reduces to $O\left(n_v^{-1}\right)$.  So the bound in (\ref{eq:skellam.error.indep}) is established.

Note that the right-hand side of (\ref{eq:stein.bnd.indep})
is analogous to the classical form of the bound for individual sums of independent indicator
random variables (e.g.,  \citep{barbour2005introduction}).  As remarked in the main text, for this
particular case of independent $Y_{ij}$, those more classical techniques could also be used
to produce the result of Theorem~\ref{thm:skellam.indep.case}.  Specifically, 
Let $T_1, T_2, \tilde T_1$, and $\tilde T_2$ be independent random variables supported on the integers.  Denote by $d_{TV}(X_1,X_2)$ the
total-variation distance between two random variables $X_1$ and $X_2$.  Then
\begin{eqnarray*}
d_{KS}\left( T_1 - T_2, \tilde T_1 - \tilde T_2\right) & \leq & d_{TV}\left(T_1 - T_2, \tilde T_1 - \tilde T_2\right) \\
 & \le & d_{TV}\left( (T_1,T_2), (\tilde T_1, \tilde T_2)\right) \\
 & \le & d_{TV}\left(T_1, \tilde T_1\right) + d_{TV}\left(T_2, \tilde T_2\right) \enskip , \\
\end{eqnarray*}
where the first inequality exploits the fact that total-variation distance provides an upper bound on Kolmogorov-Smirnov distance, 
and the second and third inequalities follow from Lemmas~3.6.3 and~3.6.2 of~\citep{durrett2010probability}, respectively.  Now define 
$$T_1 = \sum_{\left\{i,j\right\}\in E^c}Y_{ij} \quad\hbox{and}\quad T_2 =\sum_{\left\{i,j\right\}\in E}\left(1-Y_{ij}\right)\enskip ,$$
and let $\tilde T_1$ and $\tilde T_2$ be independent Poisson random variables with common mean $\lambda$.
Setting $\lambda=|E^c|\alpha=|E|\beta$, and applying to each 
of $d_{TV}\left(T_1, \tilde T_1\right)$ and $d_{TV}\left(T_2, \tilde T_2\right)$
the standard Stein bounds for Poisson approximation to sums of independent indicators (e.g., \citep[Eqn 2.6]{barbour2005introduction}),
we again obtain that
$$d_{KS}\left(D_E, \hbox{Skellam}(\lambda,\lambda)\right)\leq \frac{1}{\lambda} \left[\sum_{(i,j)\in E^c} \alpha^2+\sum_{(i,j)\in E} \beta^2\right] =\frac{|E^c| \alpha^2 + |E|\beta^2}{|E^c|\alpha}\enskip ,$$
and the rest follows.
\begin{flushright}
$\square$
\end{flushright}

\subsubsection{Proof of Theorem~\ref{thm:normal.indep.case}.}

To establish the bounds in (\ref{eq:normal.error.indep.sparse}) and (\ref{eq:normal.error.indep.dense}),
we use the following result from Stein's method for the normal distribution (e.g., \citep{barbour2005introduction}).
\begin{theorem}
Let $\xi_1,\ldots,\xi_n$ be independent random variables which have zero means and finite variances $\mathbb{E}\left[\xi_i^2\right]=\sigma_i^2$, $1\leq i \leq n$, and satisfy $\sum_{i=1}^n \sigma_i^2=1.$  If $F_n(x)$ is the cdf of $\sum_{i=1}^n \xi_i$, then, for every $\epsilon>0$,

$$\frac{1-e^{-\frac{\epsilon^2}{4}}}{40} \sum_{i=1}^n \mathbb{E}\left[\xi_i^2 I_{\left\{|\xi_i|>\epsilon\right\}}\right] - \sum_{i=1}^n \sigma_i^4 \leq \sup_{x\in \mathbb{R}} \left| F_n(x)-\Phi(x)\right| \leq 7 \sum_{i=1}^n \mathbb{E}\left[ \left| \xi_i\right|^3\right].$$
\end{theorem}

We apply this theorem, with $\xi_i=X_i/\sigma$ where $X_i$ is a term in one of the sums of $D_E$, to establish
each of our upper and lower bounds in turn.
\smallskip

{\em Upper Bounds in (\ref{eq:normal.error.indep.sparse}) and (\ref{eq:normal.error.indep.dense}):}  First, note that since
$$\sum_{i=1}^n \mathbb{E}\left[\left| \xi_i\right|^3\right] = \frac{\sum_{i=1}^n \mathbb{E}\left[\left| X_i\right|^3\right]}{\sigma^3}\enskip,$$
and
$$\mathbb{E}\left[\left| X_i\right|^3\right] = \alpha(1-\alpha) \left[(1-\alpha)^2+\alpha^2\right] \; \; \; {\rm or} \; \; \;  \beta(1-\beta) \left[(1-\beta)^2+\beta^2\right] \enskip ,$$
with $n$ understood to be either $|E^c|$ or $|E|$, it follows that
$$
\begin{aligned}
\sum_{i=1}^n \mathbb{E}\left[\left| \xi_i\right|^3\right] &=\frac{\alpha(1-\alpha) \left[(1-\alpha)^2+\alpha^2\right] |E^c| + \beta(1-\beta) \left[(1-\beta)^2+\beta^2\right] |E|}{\left(\alpha(1-\alpha)|E^c| + \beta(1-\beta)|E|\right)^{\frac{3}{2}}} \\
&\leq \max\{(1-\alpha)^2+\alpha^2,(1-\beta)^2+\beta^2\}\frac{\alpha(1-\alpha)   |E^c| + \beta(1-\beta)  |E|}{\left(\alpha(1-\alpha)|E^c| + \beta(1-\beta)|E|\right)^{\frac{3}{2}}} \\
&=\frac{\max\{(1-\alpha)^2+\alpha^2,(1-\beta)^2+\beta^2\}}{\left(\alpha(1-\alpha)|E^c| + \beta(1-\beta)|E|\right)^{\frac{1}{2}}} \\
&=\frac{\max\{(1-\alpha)^2+\alpha^2,(1-\beta)^2+\beta^2\}}{ \sqrt{2-(\alpha+\beta)}} \cdot \frac{1}{\sqrt{\alpha |E^c|}} \enskip ,
\end{aligned}
$$
where in the last equality we have used the fact that $\beta=\left(|E^c|/|E|\right) \alpha$ follows from (A2).  Finally, note that
$$(1-\alpha)^2+\alpha^2 = 1-2\alpha+2\alpha^2=1-2\alpha(1-\alpha)\leq 1$$ and the same holds for $(1-\beta)^2+\beta^2$, since $0\leq \alpha,\beta\leq 1$, so that
$$\sum_{i=1}^n \mathbb{E}\left[\left| \xi_i\right|^3\right] \leq \frac{1}{\sqrt{2-(\alpha+\beta)}} \cdot \frac{1}{\sqrt{\alpha |E^c|}}\enskip.$$
This immediately implies, after another application of $\beta=(|E^c|/|E|)\alpha$,
$$\sup_{x\in \mathbb{R}} \left| F_n(x)-\Phi(x)\right| \leq \frac{7}{\sqrt{2-(\alpha+\frac{|E^c|}{|E|}\alpha )}} \cdot \frac{1}{\sqrt{\alpha|E^c|}}\enskip.$$

Using $\alpha = \lambda / |E^c|$, and invoking the assumption of low-rate measurement error in (A3) and the definitions of sparse and dense graphs in Section~\ref{sec:initial.results}, the upper bounds in (\ref{eq:normal.error.indep.sparse}) and (\ref{eq:normal.error.indep.dense}) follow.
\smallskip

{\em Lower bound in (\ref{eq:normal.error.indep.sparse}) and (\ref{eq:normal.error.indep.dense}):}
First, note that since $\xi_i = X_i/\sigma$, $\sigma_i^2 = \alpha(1-\alpha)/\sigma^2$ or $\sigma_i^2 = \beta(1-\beta)/\sigma^2$.  Thus,
$$
\begin{aligned}
\sum_{i=1}^n \sigma_i^4  &= \frac{(\alpha(1-\alpha))^2 |E^c| + (\beta(1-\beta))^2 |E|}{(\alpha(1-\alpha)|E^c| + \beta(1-\beta)|E|)^2}\\
&= \frac{(\alpha(1-\alpha))^2 |E^c| + \left(\frac{|E^c|}{|E|} \alpha(1-\beta)\right)^2 |E|}{ (2-(\alpha+\beta))^2 (\alpha |E^c|)^2}\\
&=\frac{1}{|E^c|} \cdot \frac{1}{(2-(\alpha+\beta))^2} \cdot \left[ (1-\alpha)^2 + \frac{|E^c|}{|E|}(1-\beta)^2 \right]\\
&=\frac{1}{(2-(\alpha+\beta))^2} \left[ \frac{(1-\alpha)^2}{|E^c|} + \frac{(1-\beta)^2}{|E|}\right]\end{aligned} \enskip ,
$$
where in the second equality, we have used  $\beta=(|E^c|/|E|)\alpha$.

Next, choose $\epsilon={1}/({2\sigma})$.  Note that this is the midpoint of the intervals $$\left(  \frac{\alpha}{\sigma},\frac{1-\alpha}{\sigma}\right), \; \; \; {\rm and} \; \; \; \left(  \frac{\beta}{\sigma},\frac{1-\beta}{\sigma}\right) $$
if $\alpha,\beta<{1}/{2}$ and of the intervals
$$\left(  \frac{1-\alpha}{\sigma},\frac{\alpha}{\sigma}\right), \; \; \; {\rm and} \; \; \; \left(  \frac{1-\beta}{\sigma},\frac{\beta}{\sigma}\right)\enskip.$$
if $\alpha,\beta\geq {1}/{2}$.  In either case, these are the endpoints of the interval formed by the values of $|\xi_i|=|X_i|/\sigma$.

Due to the symmetry in these intervals about $\frac{1}{2}$ we may, without loss of generality, assume $\alpha,\beta<1/2$.
In doing so, and using $\beta=(|E^c|/|E|)\alpha$,
$$
\begin{aligned}
&\frac{1-e^{-\frac{\epsilon^2}{4}}}{40} \sum_{i=1}^n \mathbb{E}\left[\xi_i^2 I_{\left\{|\xi_i|>\epsilon\right\}}\right] \\
&=\frac{1-e^{-\frac{\epsilon^2}{4}}}{40} \cdot \frac{(1-\alpha)^2\alpha |E^c| + (1-\beta)^2 \beta|E|}{\alpha(1-\alpha)|E^c|+\beta(1-\beta)|E|}\\
&=\frac{1-e^{-\frac{\epsilon^2}{4}}}{40} \cdot \frac{(1-\alpha)^2  + (1-\beta)^2}{2-(\alpha+\beta) }\\
&=e^{-\frac{\epsilon^2}{4}} \frac{e^{\frac{\epsilon^2}{4}}-1}{40} \cdot \frac{(1-\alpha)^2  + (1-\beta)^2}{2-(\alpha+\beta) }\\
&\geq e^{-\frac{\epsilon^2}{4}} \frac{\epsilon^2}{160} \cdot \frac{(1-\alpha)^2  + (1-\beta)^2}{2-(\alpha+\beta) }\\
&= e^{-\frac{\epsilon^2}{4}} \frac{1}{640} \cdot \frac{1}{\alpha(1-\alpha)|E^c|+\beta(1-\beta)|E|} \cdot \frac{(1-\alpha)^2  + (1-\beta)^2}{2-(\alpha+\beta)}\\
&=e^{-\frac{\epsilon^2}{4}} \frac{1}{640} \cdot \frac{1}{\alpha|E^c|} \cdot \frac{(1-\alpha)^2  + (1-\beta)^2}{(2-(\alpha+\beta))^2} \enskip .
\end{aligned}
$$

Combining the two sets of expressions above, the lower bound becomes
\begin{eqnarray}
& & \frac{1-e^{-\frac{\epsilon^2}{4}}}{40} \sum_{i=1}^n \mathbb{E}\left[\xi_i^2 I_{\left\{|\xi_i|>\epsilon\right\}}\right] - \sum_{i=1}^n \sigma_i^4 \nonumber \\
&\geq & e^{-\frac{\epsilon^2}{4}} \frac{1}{640} \cdot \frac{1}{\alpha|E^c|} \cdot \frac{(1-\alpha)^2  + (1-\beta)^2}{(2-(\alpha+\beta))^2} - \frac{1}{(2-(\alpha+\beta))^2} \left[ \frac{(1-\alpha)^2}{|E^c|} + \frac{(1-\beta)^2}{|E|}\right] \nonumber \\
&=& \frac{1}{(2-(\alpha+\beta))^2} \left[ (1-\alpha)^2 \left( \frac{\exp\left(-\frac{1}{16} \cdot \frac{1}{\alpha|E^c|}\cdot \frac{1}{2-(\alpha+\beta)}\right)}{640} \cdot \frac{1}{\alpha|E^c|} - \frac{1}{|E^c|}\right) \right. \nonumber \\
& &\hspace{1in}\left. +(1-\beta)^2 \left( \frac{\exp\left(-\frac{1}{16} \cdot \frac{1}{\alpha|E^c|}\cdot \frac{1}{2-(\alpha+\beta)}\right)}{640} \cdot \frac{1}{\alpha|E^c|} - \frac{1}{|E|}\right)\right] \enskip.
\label{eq:ugly.lower.bnd}
\end{eqnarray}
But for sufficiently large $n_v$, the exponential term in (\ref{eq:ugly.lower.bnd}) behaves like
$\exp\left[-1/(16\lambda)\right] \approx 1 - (1/16\lambda)$.  Substituting accordingly and simplifying to ignore the various terms tending to a constant in large $n_v$, the expression in (\ref{eq:ugly.lower.bnd}) can be seen to behave asymptotically like
\begin{eqnarray}
\frac{1}{4}\left[ \left(\frac{1}{640\lambda} - \frac{1}{|E^c|}\right) +
                        \left(\frac{1}{640\lambda} - \frac{1}{|E|}\right)\right] \enskip .
\label{eq:pretty.lower.bnd}
\end{eqnarray}

Again, by the assumption of low-rate measurement error in (A3) and the definitions of sparse and dense graphs given in Section~\ref{sec:initial.results}, appropriate substitution of the values for $\lambda$, $|E|$, and $|E^c|$ yield the lower bounds in (\ref{eq:normal.error.indep.sparse}) and (\ref{eq:normal.error.indep.dense}).  This completes the proof of Theorem~\ref{thm:normal.indep.case}.
\begin{flushright}
$\square$
\end{flushright}

\subsection{Proof of Theorem~\ref{thm:bi.monotone.coupling}.}
\label{sec:bi.monotone.coupling}

The proof follows by rewriting each of the two sums bracketed in (\ref{Stein KS Bound}), and then aggregating
terms.  Under condition $(i)$ of the theorem,
\begin{eqnarray*}
\left|U-U_k^{(L)} \right| & = & \left|L_k + \sum_{j\ne k} L_j - \sum_{\ell} M_\ell
                                                   - \left( \sum_{j\ne k} \tilde{L}^{L_k}_j - \sum_{\ell} \tilde{M}^{L_k}_\ell\right)\right| \cr
& = & \sum_{j\ne k} \tilde{L}^{L_k}_j - \sum_{\ell} \tilde{M}^{L_k}_\ell - 
         \left(\sum_{j\ne k} L_j - \sum_{\ell} M_\ell\right) - L_k \cr
& = & U_k^{(L)} - U \enskip .
\end{eqnarray*}
Similarly, under condition $(ii)$ of the theorem, $\left| U-U_k^{(M)} \right| = U-U_k^{(M)}$.

In the absence of having to deal directly with the absolute values, we find that
$$\sum_{k=1}^n p_k E\left|U-U_k^{(L)} \right| =
\sum_{k=1}^n \sum_{k\ne j} E\left[L_kL_j\right] - \sum_{k,\ell} E\left[L_k M_\ell\right] 
                - E\left[U\right] \lambda_1$$
and
$$\sum_{k=1}^m q_k E\left| U-U_k^{(M)} \right| = 
\sum_{k=1}^m \sum_{\ell\ne k} E\left[M_k M_\ell\right] - \sum_{k,\ell} E\left[L_k M_\ell\right]
                 + E\left[U\right]\lambda_2 \enskip .$$
As a result, the bracketed term in (\ref{Stein KS Bound}) takes the form
\begin{eqnarray*}
& & \sum_{k=1}^n p_k E\left|U-U_k^{(L)} \right|  + \sum_{k=1}^m q_k E\left| U-U_k^{(M)} \right| \cr
& = & \sum_{k=1}^n \sum_{k\ne j} E\left[L_kL_j\right] + \sum_{k=1}^m \sum_{\ell\ne k} E\left[M_k M_\ell\right] 
         - 2 \sum_{k,\ell} E\left[L_k M_\ell\right] - \left(E\left[U\right]\right)^2 \cr
& = & E\left[U^2\right] - \left(E\left[ U\right]\right)^2 
         - \sum_{k=1}^n E\left[L_k^2\right] - \sum_{\ell=1}^m E\left[M^2_\ell\right] \cr
& = & \hbox{Var}(U) - \sum_{k=1}^n E\left[L_k\right] - \sum_{\ell=1}^m E\left[M_\ell\right] \cr
& = & \hbox{Var}(U) - \left(\lambda_1+\lambda_2\right) \cr
& = & \hbox{Var}(U) - \hbox{Var}(W)\enskip .
\end{eqnarray*}
\begin{flushright}
$\square$
\end{flushright}



\bibliographystyle{plainnat}
\bibliography{bibliog}

\end{document}